\documentclass[10pt, reqno]{amsart}

\usepackage{amsmath, amssymb, amsthm, amsaddr}
\usepackage{geometry}
\usepackage{tikz}
\usepackage{booktabs}
\usepackage{comment}
\usepackage{longtable}

\IfFileExists{fontawesome5.sty}{\usepackage{fontawesome5}}
                               {}
\usepackage{xcolor}
\usepackage[colorlinks=true, linkcolor=blue, urlcolor=blue,
            citecolor=blue]{hyperref}

\newtheorem{thm}{Theorem}[section]
\newtheorem{lemma}[thm]{Lemma}
\newtheorem{cor}[thm]{Corollary}
\newtheorem{conjecture}[thm]{Conjecture}
\newtheorem{ques}[thm]{Question}

\theoremstyle{definition}
\newtheorem{defn}[thm]{Definition}
\newtheorem{example}[thm]{Example}
\theoremstyle{definition}
\newtheorem{construction}[thm]{Construction}

\theoremstyle{remark}
\newtheorem{remark}[thm]{Remark}
\newtheorem{observation}[thm]{Observation}

\numberwithin{equation}{section}

\newcommand{\Z}{\mathbb{Z}}
\newcommand{\V}{\mathcal{V}}
\newcommand{\calP}{\mathcal{P}}
\newcommand{\Lcal}{\mathcal{L}}
\newcommand{\Hcal}{\mathcal{H}}
\newcommand{\Rcal}{\mathcal{R}}
\newcommand{\Fcal}{\mathcal{F}}

\newcommand{\rthmname}{Theorem}
\newtheorem*{rthm}{\rthmname}
\newenvironment{restatethm}[1]{\renewcommand{\rthmname}{#1}
\begin{rthm}
}{\end{rthm}}

\tikzset{
  vtx/.style={circle, draw, fill=white, minimum size=18pt, inner sep=0pt,
              font=\small},
  ed/.style={black!70, thick},
  blk/.style={font=\small\itshape, gray}
}

\newcommand{\tworows}[2]{%
  \foreach \lbl [count=\i from 0] in {#1} {\node[vtx] (t\i) at (\i,1) {$\lbl$};}%
  \foreach \lbl [count=\i from 0] in {#2} {\node[vtx] (b\i) at (\i,0) {$\lbl$};}%
}

\newcommand{\blocklabels}[3]{%
  \node[blk] at (1,-0.6) {#1};
  \node[blk] at (4,-0.6) {#2};
  \node[blk] at (7,-0.6) {#3};
  \draw[dashed, gray!40] (2.5,-1) -- (2.5,1.5);
  \draw[dashed, gray!40] (5.5,-1) -- (5.5,1.5);
}

\title{Perfect Partitions of Regular Bipartite Graphs}

\date{\today}

\author{Melanie J. Ferreri, Archana L. Parameswaran, Gexin Yu}
\address{Department of Mathematics, William \& Mary, Williamsburg, VA 23185
}

\email {mjferreri@wm.edu, aparameswaran9@gatech.edu, gyu@wm.edu}

\begin{document}
\maketitle

\begin{abstract}
Let $I_{n,n}$ be a fixed perfect matching of the complete bipartite graph
$K_{n,n}$, and let $L_{n,1}=K_{n,n}-I_{n,n}$. A \emph{perfect set} of $L_{n,1}$
is a $1$-factorization of $L_{n,1}$, that is, a set of $n-1$ pairwise
edge-disjoint perfect matchings whose union is $E(L_{n,1})$; a \emph{perfect
partition} of $L_{n,1}$ is a partition of the set of all perfect matchings of
$L_{n,1}$ into perfect sets. Brualdi, Chiang, and Li (2005) asked whether
$L_{n,1}$ admits a perfect partition. Identifying both sides of the bipartition
with $[n]=\{1,2,\ldots ,n\}$, the perfect matchings of $L_{n,1}$ are exactly the
derangements of $[n]$, the perfect sets correspond bijectively to the reduced
Latin squares of order $n$, and the perfect partitions are exactly the perfect
matchings of the $(n-1)$-uniform hypergraph whose vertices are the derangements
of $[n]$ and whose edges are the perfect sets. We first give, for every $n$, a
constructive lower bound on the maximum number of pairwise disjoint perfect sets
of $L_{n,1}$. We then show that $L_{5,1}$ has exactly $42$ distinct perfect partitions and
that $L_{6,1}$ has at least $82$ distinct perfect partitions; these are the first examples demonstrating that
perfect partitions need not be unique.
Finally, we exhibit two distinct perfect partitions of $L_{7,1}$, resolving the
 $n=7$ open problem raised by Brualdi et al.
\end{abstract}

\section{Introduction}\label{sec:intro}

Let $n\ge 2$. We take the labeled complete bipartite graph $K_{n,n}$ to have
bipartition $(V,V')$ with $V=V'=[n]:=\{1,2,\ldots ,n\}$, we write $ij'$ for the
edge joining $i\in V$ to $j'\in V'$, and we let $I_{n,n}=\{\,ii' : i\in[n]\,\}$
be the perfect matching joining equal labels. Define
\[
  L_{n,1} \;=\; K_{n,n}-I_{n,n},
\]
so that $L_{n,1}$ is $(n-1)$-regular and $|E(L_{n,1})|=n(n-1)$.

\begin{defn}\label{def:perfect}
A \emph{perfect set} of $L_{n,1}$ is a $1$-factorization of $L_{n,1}$, that is, a
set $\{M_1,\ldots ,M_{n-1}\}$ of pairwise edge-disjoint perfect matchings of
$L_{n,1}$ satisfying $\bigcup_{i}M_i=E(L_{n,1})$. Two perfect sets are
\emph{disjoint} if they contain no common perfect matching. A \emph{perfect
partition} of $L_{n,1}$ is a partition of the set of all perfect matchings of
$L_{n,1}$ into perfect sets.
\end{defn}

Every $1$-factorization of $L_{n,1}$ has exactly $(n-1)$ 1-factors, since $L_{n,1}$
is $(n-1)$-regular. Distinct perfect sets always share edges, so disjointness in
Definition~\ref{def:perfect} means disjointness as sets of perfect matchings, not disjointness as sets of edges.

\begin{ques}[Brualdi, Chiang, and Li \cite{BCL05}]\label{ques:BCL}
Does $L_{n,1}$ have a perfect partition for every $n \ge 2$?
\end{ques}

\subsection{Equivalent formulations}\label{subsec:formulations}

Question~\ref{ques:BCL} can be posed in a variety of different ways. We describe four formulations below.

\medskip
\noindent\textbf{Permutation matrices.}
The biadjacency matrix of $L_{n,1}$ is $L_n=J_n-I_n$, where $J_n$ is the
all-ones matrix and $I_n$ is the identity matrix.
The permutation matrices $P$ of order $n$ with $P\le L_n$ entrywise are exactly
the permutation matrices with zero diagonal, and these correspond to the perfect
matchings of $L_{n,1}$. A perfect set is therefore a decomposition
$L_n=P_1+P_2+\cdots +P_{n-1}$ into permutation matrices, and
Question~\ref{ques:BCL} asks whether the set of all zero-diagonal permutation
matrices of order $n$ can be partitioned into such decompositions. Brualdi et al. first studied the problem of perfect partitions in this context \cite{BCL05}.

\medskip

\noindent\textbf{Derangements.}
Let $S_n$ denote the group of permutations of $[n]$. Under the identification $V=V'=[n]$, a perfect matching $M$ of $K_{n,n}$ is the
graph of a permutation $\sigma\in S_n$, namely
$M=\{\,i\,\sigma(i)' : i\in[n]\,\}$, and $M\subseteq E(L_{n,1})$ if and only if
$\sigma(i)\ne i$ for all $i$. Thus the perfect matchings of $L_{n,1}$ are exactly
the derangements of $[n]$, the set of which we denote by $D_n$, and let $d_n$ denote the size, $d_n = |D_n|$. For derangements 
$\sigma,\tau\in S_n$ we observe when $\sigma(i)\ne\tau(i)$ for every
$i\in[n]$, and write $\sigma\perp\tau$ if that is the case. Equivalently we look at when $\tau^{-1}\sigma$ is a derangement, or when $\sigma$ and $\tau$ are \emph{mutually deranged}. Two perfect
matchings of $L_{n,1}$ are edge-disjoint precisely when the corresponding
derangements are mutually deranged. Consequently a set of $n-1$ pairwise
mutually deranged derangements is automatically a perfect set, since its members
contribute $n(n-1)=|E(L_{n,1})|$ distinct edges. Hence
\[
  L_{n,1}\ \text{has a perfect partition}
  \iff
  \begin{minipage}{0.45\textwidth}
  the derangements of $[n]$ can be partitioned into sets of $n-1$ pairwise
  mutually deranged derangements.
  \end{minipage}
\]
It is known that the number of derangements satisfies the recurrence relation $d_n=(n-1)(d_{n-1}+d_{n-2})$, with $d_0 = 1$ and $d_1 = 0$.
From this
we see $(n-1)\mid d_n$, so there is no divisibility obstruction. As a result, every perfect partition of $L_{n,1}$
consists of exactly $d_{n-1}+d_{n-2}$ perfect sets.

\medskip
\noindent\textbf{Reduced Latin squares.}
A \emph{reduced Latin square} of order $n$ is an $n\times n$ array on the symbol
set $[n]$ in which every row and every column contains each symbol exactly once,
and whose first row and first column are $1,2,\ldots ,n$ in that order. 
As Lemma~\ref{lem:latin} shows, deleting the first row of such a square leaves a perfect set, and the set of all reduced Latin squares of order $n$ provides the list of all possible perfect sets of $L_{n,1}$.

\begin{lemma}\label{lem:latin}
For a reduced Latin square $R$ of order $n$ and $2\le r\le n$, put
$\sigma_r(j)=R(r,j)$. Then $R\mapsto\{\sigma_2,\dots ,\sigma_n\}$ is a bijection
from the reduced Latin squares of order $n$ onto the perfect sets of $L_{n,1}$.
\end{lemma}

\begin{proof}
Let $R$ be a reduced Latin square and $2\le r\le n$. Each row of $R$ is a
permutation of $[n]$, so $\sigma_r\in S_n$. Since $R(1,j)=j$ and column $j$
contains each symbol exactly once, $\sigma_r(j)\ne j$ for every $j$, so
$\sigma_r$ is a derangement; the same column condition gives
$\sigma_r(j)\ne\sigma_{r'}(j)$ for $r\ne r'$, so
$\sigma_r\perp\sigma_{r'}$. Thus $\{\sigma_2,\dots ,\sigma_n\}$ is a set of
$n-1$ pairwise mutually deranged derangements, hence a perfect set.

Conversely let $\{\tau_1,\dots ,\tau_{n-1}\}$ be a perfect set. The values
$\tau_i(1)$ are pairwise distinct and lie in $[n]\setminus\{1\}$, so they exhaust
it. Build an $n\times n$ array whose first row is $1,2,\dots ,n$ and whose row
$\tau_i(1)$ is $\tau_i(1),\tau_i(2),\dots ,\tau_i(n)$. Every row is a
permutation of $[n]$. Column $j$ contains $j$ once, from row $1$, together with
the $n-1$ values $\tau_i(j)$, which are pairwise distinct and all different from
$j$ and therefore exhaust $[n]\setminus\{j\}$; so every column contains each
symbol once. The first column reads $1,2,\dots ,n$ by construction. Hence the
array is a reduced Latin square, and the two maps are mutually inverse.
\end{proof}

We call two reduced Latin squares \emph{disjoint} if the corresponding perfect sets
are disjoint. Then
\[
  L_{n,1}\ \text{has a perfect partition}
  \iff
  \text{There exist $\tfrac{d_n}{n-1}$ pairwise disjoint reduced Latin squares.}
\]

\medskip
\noindent\textbf{Hypergraph matchings.}
Let $\Hcal_n=(V,E)$ be the hypergraph with $V(\Hcal_n)=D_n$, the set of derangements of $[n]$. We define the edges of $\Hcal_n$ as follows
\[
  E(\Hcal_n)=\bigl\{\{v_1,\dots ,v_{n-1}\}: v_i\in D_n
  \text{ pairwise mutually deranged}\bigr\},
\]

such that the hyperedges are exactly the perfect sets. Then $\Hcal_n$ is
$(n-1)$-uniform with $d_n$ vertices, a matching in $\Hcal_n$ is a family of
pairwise disjoint perfect sets, and a perfect partition of $L_{n,1}$ is exactly a
perfect matching of $\Hcal_n$.
Writing $\nu(\Hcal_n)$ for the matching number,
\[
  L_{n,1}\ \text{has a perfect partition}
  \iff
  \nu(\Hcal_n) = d_{n-1}+d_{n-2}.
\]

\subsection{Known results and our contribution}\label{subsec:results}

Brualdi, Chiang, and Li \cite{BCL05} answered Question~\ref{ques:BCL} for
$n\le 5$, and Li, Soosiah, and Yu \cite{LSY12} found a perfect partition for the $L_{6,1}$ case. The existence of a perfect partition of
$L_{n,1}$ for $n\ge 7$ was open. We solve the open problem of a perfect partition for $n=7$ (Theorem~\ref{thm:small_cases}), and $n \geq 8$ remains open. The difficulty with finding perfect partitions beyond small cases stems from a lack of algebraic
structure we can exploit from our perfect matchings, and the number of
candidate perfect sets (equivalently, the number of reduced Latin squares of order $n$) grows very rapidly, so exhaustive search becomes infeasible quickly.

Throughout, we let $\lvert\calP(L_{n,1})\rvert$ denote the
maximum number of pairwise disjoint perfect sets of $L_{n,1}$.
Theorem~\ref{thm:general_lower_bound} gives a constructive lower bound on
$\lvert\calP(L_{n,1})\rvert$ for all $n$. Theorem~\ref{thm:small_cases} records
the counts of perfect partitions obtained for $n\in\{5,6,7\}$; these show in
particular that perfect partitions are not unique, and answer the question of
\cite{BCL05} for $n=7$.  

\begin{thm}\label{thm:general_lower_bound}
Let $n\ge 2$, let $p$ be the largest prime factor of $n$, set $k=n/p$, and let
$|\mathcal P(L_{p,1})|$ denote the maximum number of pairwise disjoint perfect sets of
$L_{p,1}$. Then
\[
  \lvert \calP(L_{n,1})\rvert
  \;\ge\;
  \frac{n!}{(p!)^{k}\, k!}\;|\mathcal P(L_{p,1})|\bigl((p-2)!\bigr)^{k-1} (p-1)^{\,k-1}.
\]
\end{thm}

\begin{thm}\label{thm:small_cases}
$L_{5,1}$ has exactly $42$ perfect partitions, $L_{6,1}$ has at least $82$
perfect partitions, and $L_{7,1}$ has at least $2$ perfect partitions.
\end{thm}

Section~\ref{sec:blocks} introduces the block graphs $\Lcal_{p,1;k}$ and
produces many pairwise disjoint perfect sets inside each of them.
Section~\ref{sec:layers} adds the crossing edges, in the form of layer graphs
$Q^{\kappa}_{\V}$, and combines the two constructions to prove
Theorem~\ref{thm:general_lower_bound}.
Section~\ref{sec:small} describes the graph model used in our computational
search to prove our results for $n=5,6,7$ and proving Theorem~\ref{thm:small_cases}.
Section~\ref{sec:future} discusses open problems, including what is known about
$\Hcal_n$ asymptotically.

\section{Disjoint perfect sets of the block graphs \texorpdfstring{$\Lcal_{p,1;k}$}{L(p,1;k)}}
\label{sec:blocks}

Throughout this section $p$ is a prime, $k\ge1$, and $n=pk$. We partition $[n]$
into $k$ blocks of size $p$ and place a copy of $L_{p,1}$ on each block; the
resulting spanning subgraphs of $L_{n,1}$ are the objects of study.

\subsection{Vertex configurations}\label{subsec:configs}

\begin{defn}\label{def:config}
A \emph{vertex configuration} is an ordered partition
$\V=(V_1,\dots ,V_k)$ of $[n]$ into $k$ blocks with $|V_i|=p$ for every $i$; we
write $V_i'=\{\,x' : x\in V_i\,\}$ for the corresponding set of primed vertices.
The graph $L_{\V}$ is the spanning subgraph of $L_{n,1}$ with edge set
\[
  E(L_{\V})\;=\;\bigl\{\,xy' \;:\; x,y\in V_i \text{ for some } i,\ x\ne y \,\bigr\},
\]
that is, an edge copy of $L_{p,1}$ on each vertex pair $V_i/V_i'$. We write $\Lcal_{p,1;k}$
for the family of all such graphs.
\end{defn}

\begin{example}\label{ex:LV}
Let $L_{\V}\in \Lcal_{3,1;3}$ with
$\V= (V_1 = \{1, 5, 6\},\; V_2 = \{2, 8, 9\},\; V_3 = \{3, 4, 7\})$. Applying
the $L_{3,1}$ edge set to each block gives
\[
L_\V \;=\;
\begin{tikzpicture}[baseline=(current bounding box.center)]
  \tworows{1,5,6,2,8,9,3,4,7}{1',5',6',2',8',9',3',4',7'}
  \draw[ed] (t0)--(b1) (t0)--(b2) (t1)--(b0) (t1)--(b2) (t2)--(b0) (t2)--(b1);
  \draw[ed] (t3)--(b4) (t3)--(b5) (t4)--(b3) (t4)--(b5) (t5)--(b3) (t5)--(b4);
  \draw[ed] (t6)--(b7) (t6)--(b8) (t7)--(b6) (t7)--(b8) (t8)--(b6) (t8)--(b7);
  \blocklabels{$V_1/V_1'$}{$V_2/V_2'$}{$V_3/V_3'$}
\end{tikzpicture}
\]
\end{example}

\begin{remark}\label{rem:componentwise}
We use componentwise notation for these graphs: $L_{\V}\in\Lcal_{p,1;3}$ is
written $L_{\V}=(L_{p,1},L_{p,1},L_{p,1})$, the $i$-th component recording the
edges on the pair $V_i/V_i'$. A perfect matching of $L_{\V}$ is likewise written
as a $k$-tuple of perfect matchings of the individual copies of $L_{p,1}$. Additionally, although we draw dashed separations for clarity in the image, $L_\V$ is a single graph.
\end{remark}

\subsection{Perfect sets of blocks}\label{subsec:pcycles}

The block input to the construction is a supply of pairwise disjoint perfect
sets of $L_{p,1}$. In particular, the $p$-cycles provide one such supply because they have a nice group structure.

\begin{lemma}\label{lem:p-cycles}
Let $p$ be prime. The $(p-1)!$ $p$-cycles in $S_p$ can be partitioned into
$(p-2)!$ disjoint sets, each containing $p-1$ elements whose permutation
matrices sum to $L_p=J_p-I_p$. Thus the set of all $p$-cycle perfect matchings
of $L_{p,1}$ can be partitioned into perfect sets.
\end{lemma}

\begin{proof}
Let $\sigma \in S_p$ be a $p$-cycle. Since $p$ is prime, the cyclic group
$\langle \sigma \rangle$ has order $p$ and every non-identity element of it is a
$p$-cycle. By Lagrange's theorem, distinct subgroups of order $p$ intersect
trivially. Therefore the sets $\langle \sigma \rangle \setminus \{e\}$ partition
the $(p-1)!$ $p$-cycles into $\frac{(p-1)!}{p-1} = (p-2)!$ disjoint sets. The
non-identity elements of each subgroup act regularly on $[p]$, so their
permutation matrices sum to $J_p - I_p = L_p$.
\end{proof}

\begin{defn}\label{def:PCp}
Let $\calP(C_p)$ denote the set of $(p-2)!$ perfect sets of $L_{p,1}$ produced by
Lemma~\ref{lem:p-cycles}.

\end{defn}

Define
\[\mathfrak{P}(L_{p,1}): = \{\mathcal{P}=\{P_1, P_2, \ldots, P_j\}: \mathcal{P} \text{ is a set of pairwise disjoint perfect sets on } L_{p,1}\}.\]

We are interested in the elements of largest size in $\mathfrak{P}(L_{p,1})$, that is the largest set of pairwise disjoint perfect sets we have for each $p$. Although there could exist multiple sets of pairwise disjoint sets of maximum size, we select and fix an arbitrary one, which we denote by $\mathcal P(L_{p,1}) \in \mathfrak{P}(L_{p,1})$. Then, note
\[|\mathcal{P}(L_{p,1})|=\displaystyle\max_{\mathcal{P} \in \mathfrak{P}(L_{p,1})} | \mathcal{P}|\]

Lemma~\ref{lem:p-cycles} gives $|\mathcal P(L_{p,1})|\ge(p-2)!$, while
$|\mathcal P(L_{p,1})|\le d_p/(p-1)$ .

\begin{remark}\label{rem:ordering}
For every perfect set in $\calP(C_p)$ and in $\mathcal P(L_{p,1})$, we fix an
arbitrary ordering of its perfect matchings and write the perfect set as an ordered tuple
$P=(M^{1},\dots ,M^{p-1})$. This ordering stays the same throughout, so that
the rotations of Definition~\ref{def:rotation} are well defined.
\end{remark}

\subsection{Rotations}\label{subsec:rotations}

\begin{defn}\label{def:rotation}
Let $L_{\V}\in\Lcal_{p,1;k}$, and choose $P_1\in\mathcal P(L_{p,1})$ and
$P_i\in\calP(C_p)$ for $2\le i\le k$, placing the edge set of $P_i$ on the pair
$V_i/V_i'$. Using the ordering of Remark~\ref{rem:ordering}, write
$P_i=(M_i^{1},\dots ,M_i^{p-1})$, so that $M_i^{j}$ is the $j$-th member of the
perfect set on block $V_i$. Let $\rho=(1\ 2\ \cdots\ p-1)$ be the cyclic
permutation of $\{1,\dots ,p-1\}$ of order $p-1$. For a \emph{rotation vector}
$\mathbf r=(r_1,\dots ,r_k)$ with $r_1=0$ and
$r_2,\dots ,r_k\in\{0,\dots ,p-2\}$, define
\begin{equation}\label{eq:rotation}
  R^{\V}_{\mathbf r}(P_1,\dots ,P_k)
  \;=\;
  \Bigl\{\,\bigl(M_1^{\ell},\,M_2^{\rho^{r_2}(\ell)},\,\dots ,\,
                  M_k^{\rho^{r_k}(\ell)}\bigr)
  \;:\; \ell=1,\dots ,p-1 \,\Bigr\},
\end{equation}
a set of $p-1$ perfect matchings of $L_{\V}$ in the componentwise notation of
Remark~\ref{rem:componentwise}. The first coordinate is never rotated, since fixing
$r_1=0$ fixes the ordering of the first block's perfect set, so we reduce $\rho^{r_1}(\ell)$ to simply $\ell$.
\end{defn}

\begin{example}\label{ex:rotation}
From Example~\ref{ex:LV}, let $L_{\V}\in\Lcal_{3,1;3}$ with
$V_1=\{1,5,6\}$, $V_2=\{2,8,9\}$, $V_3=\{3,4,7\}$, so $p=3$, $k=3$, and $p-1=2$.
Here $\rho=(1\ 2)$, so $\rho^0=\mathrm{id}$ and $\rho^1=(1\ 2)$, and $r_1=0$ while $r_2,r_3\in\{0,1\}$.

Take $P_1,P_2,P_3$ to be the perfect set $\{M^1,M^2\}$ of $3$-cycles on each block, where
$M_i^1$ is the cycle through the elements of $V_i$ in the order listed above. The following labeling is now fixed for our construction:
\begin{alignat*}{2}
  M_1^1&=(1\ 5\ 6), &\qquad M_1^2&=(1\ 6\ 5),\\
  M_2^1&=(2\ 8\ 9), &\qquad M_2^2&=(2\ 9\ 8),\\
  M_3^1&=(3\ 4\ 7), &\qquad M_3^2&=(3\ 7\ 4).
\end{alignat*}
Form the rotation vector $\mathbf r=(0,0,0)$, where $r_1=r_2=r_3=0$, so $\rho^{r_2}= \rho^{r_3}=\mathrm{id}$. We have no rotation here, thus we plainly apply these edge sets to form the following perfect matching of $ L_\V$:
\[
(M_1^1,M_2^1,M_3^1) \;=\;
\begin{tikzpicture}[baseline=(current bounding box.center)]
  \tworows{1,5,6,2,8,9,3,4,7}{1',5',6',2',8',9',3',4',7'}
  \draw[ed] (t0)--(b1) (t1)--(b2) (t2)--(b0);
  \draw[ed] (t3)--(b4) (t4)--(b5) (t5)--(b3);
  \draw[ed] (t6)--(b7) (t7)--(b8) (t8)--(b6);
  \blocklabels{$V_1/V_1'$}{$V_2/V_2'$}{$V_3/V_3'$}
\end{tikzpicture}
\]
\[
(M_1^2,M_2^2,M_3^2) \;=\;
\begin{tikzpicture}[baseline=(current bounding box.center)]
  \tworows{1,5,6,2,8,9,3,4,7}{1',5',6',2',8',9',3',4',7'}
  \draw[ed] (t0)--(b2) (t2)--(b1) (t1)--(b0);
  \draw[ed] (t3)--(b5) (t5)--(b4) (t4)--(b3);
  \draw[ed] (t6)--(b8) (t8)--(b7) (t7)--(b6);
  \blocklabels{$V_1/V_1'$}{$V_2/V_2'$}{$V_3/V_3'$}
\end{tikzpicture}
\]
Now we choose a different rotation vector of $\mathbf r=(0,1,0)$, where $\rho^{r_2}=(1\ 2)$ and $\rho^{r_1}=\rho^{r_3}=\mathrm{id}$. This forms the perfect set
$R^{\V}_{\mathbf r}(P_1,P_2,P_3)=\{(M_1^1,M_2^2,M_3^1),\,(M_1^2,M_2^1,M_3^2)\}$. That is,
the second block is rotated, while the first and
third are not. We created this new pairwise disjoint perfect set of $L_\V$: 
\[
(M_1^1,M_2^2,M_3^1) \;=\;
\begin{tikzpicture}[baseline=(current bounding box.center)]
  \tworows{1,5,6,2,8,9,3,4,7}{1',5',6',2',8',9',3',4',7'}
  \draw[ed] (t0)--(b1) (t1)--(b2) (t2)--(b0);
  \draw[ed] (t3)--(b5) (t5)--(b4) (t4)--(b3);
  \draw[ed] (t6)--(b7) (t7)--(b8) (t8)--(b6);
  \blocklabels{$V_1/V_1'$}{$V_2/V_2'$}{$V_3/V_3'$}
\end{tikzpicture}
\]
\[
(M_1^2,M_2^1,M_3^2) \;=\;
\begin{tikzpicture}[baseline=(current bounding box.center)]
  \tworows{1,5,6,2,8,9,3,4,7}{1',5',6',2',8',9',3',4',7'}
  \draw[ed] (t0)--(b2) (t2)--(b1) (t1)--(b0);
  \draw[ed] (t3)--(b4) (t4)--(b5) (t5)--(b3);
  \draw[ed] (t6)--(b8) (t8)--(b7) (t7)--(b6);
  \blocklabels{$V_1/V_1'$}{$V_2/V_2'$}{$V_3/V_3'$}
\end{tikzpicture}
\]
\end{example}

\begin{construction}\label{construction}\leavevmode

\begin{enumerate}
  \item Factor $n=\prod_{j=1}^{s}q_j^{e_j}$ into primes. Let $p$ be the largest
        prime factor and set $k=n/p$. Consider only $\Lcal_{p,1;k}$.
  \item Fix a partition $\V=\{V_1,\dots ,V_k\}$ of $[n]$ into $k$ blocks of size
      $p$, indexed so that $\min \{V_1\}<\cdots<\min \{V_k\}$; in particular $1\in V_1$.
      Each $V_i$ is identified with $[p]$ by the order-preserving bijection.
  \item Choose perfect sets $P_1,\dots ,P_k$ with $P_1$ one of the $|\mathcal P(L_{p,1})|$
        pairwise disjoint perfect sets of $L_{p,1}$ on $V_1/V_1'$, and
        $P_i\in\calP(C_p)$ for $i\ge2$.
  \item Choose a rotation vector $\mathbf r$ with $r_1=0$ and
        $r_2,\dots ,r_k\in\{0,\dots ,p-2\}$, and form
        $R^{\V}_{\mathbf r}(P_1,\dots ,P_k)$.
\end{enumerate}

\end{construction}

\subsection{The construction produces disjoint perfect sets} \label{subsec:disjointness}

\begin{lemma}\label{lem:rotation-is-perfect}
$R^{\V}_{\mathbf r}(P_1,\dots ,P_k)$ from Construction \ref{construction} is a perfect set of $L_{\V}$.
\end{lemma}

\begin{proof}
Each tuple $\bigl(M_1^{\ell},M_2^{\rho^{r_2}(\ell)},\dots ,
M_k^{\rho^{r_k}(\ell)}\bigr)$ matches the vertices of $V_i$ bijectively to those
of $V_i'$ for every $i$, and the blocks partition $[n]$, so each tuple is a
perfect matching of $L_{\V}$.

Since the blocks are vertex-disjoint, the edge union of the $p-1$ tuples may be
taken componentwise:
\[
\bigcup_{\ell=1}^{p-1}
\bigl(M_1^{\ell},\dots ,M_k^{\rho^{r_k}(\ell)}\bigr)
=
\Bigl(\ \bigcup_{\ell=1}^{p-1}M_1^{\ell},\ \
        \bigcup_{\ell=1}^{p-1}M_2^{\rho^{r_2}(\ell)},\ \ \dots ,\ \
        \bigcup_{\ell=1}^{p-1}M_k^{\rho^{r_k}(\ell)}\ \Bigr).
\]
Every power of $\rho$ is a bijection of $\{1,\dots ,p-1\}$, so
$\bigcup_{\ell}M_i^{\rho^{r_i}(\ell)}=\bigcup_{\ell}M_i^{\ell}$, which is the
edge set of the copy of $L_{p,1}$ on $V_i/V_i'$ because $P_i$ is a perfect set
there. Hence the union of the $p-1$ tuples is $E(L_\V)$.

They are also pairwise edge-disjoint. Let $\ell\ne\ell'$ and fix $i$. Since
$\rho^{r_i}$ is injective, $\rho^{r_i}(\ell)\ne\rho^{r_i}(\ell')$, so the $i$-th
components $M_i^{\rho^{r_i}(\ell)}$ and $M_i^{\rho^{r_i}(\ell')}$ are distinct
members of the perfect set $P_i$ and hence share no edge. This holds for every
$i$, and every edge of $L_\V$ lies in some block, so the two tuples share no
edge.
\end{proof}

\begin{lemma}\label{lem:diff-P}
Fix a vertex configuration $\V$ for $L_{\V}\in\Lcal_{p,1;k}$ and from Construction \ref{construction} create $R^{\V}_{\mathbf r}(P_1,\dots ,P_k)$ and $R^{\V}_{\mathbf r'}(P_1',\dots ,P_k')$  where $P_i\ne P_i'$ for some
$i$. Then
$R^{\V}_{\mathbf r}(P_1,\dots ,P_k)$ and $R^{\V}_{\mathbf r'}(P_1',\dots ,P_k')$
are disjoint perfect sets of $L_{\V}$.
\end{lemma}

\begin{proof}
Suppose for contradiction that the two sets share a perfect matching, say
\[
  \bigl(M_1^{\ell},M_2^{\rho^{r_2}(\ell)},\dots ,M_k^{\rho^{r_k}(\ell)}\bigr)
  =
  \bigl(N_1^{\ell'},N_2^{\rho^{r_2'}(\ell')},\dots ,N_k^{\rho^{r_k'}(\ell')}\bigr),
\]
where $M_i^{\bullet}$ denotes a member of $P_i$ and $N_i^{\bullet}$ a member of
$P_i'$. Since $L_{\V}$ is a disjoint union over the blocks, equality holds
componentwise, so in particular the $i$-th components agree. But $P_i$ and
$P_i'$ are distinct members of $\mathcal{P}(L_{p,1})$ (or $\mathcal{P}(C_p)$), hence disjoint, so they share no
perfect matching, which is a contradiction.
\end{proof}

\begin{lemma}\label{lem:diff-r}
Fix a vertex configuration $\V$ and from Construction \ref{construction} create $R^{\V}_{\mathbf r}(P_1,\dots ,P_k)$ and
$R^{\V}_{\mathbf r'}(P_1,\dots ,P_k)$ 
using rotation vectors $\mathbf r$ and $\mathbf r'$ where $r_i\ne r_i'$ for some $i$. Then $R^{\V}_{\mathbf r}(P_1,\dots ,P_k)$ and
$R^{\V}_{\mathbf r'}(P_1,\dots ,P_k)$ are disjoint perfect sets of $L_{\V}$.
\end{lemma}

\begin{proof}
Suppose for contradiction that for some $\ell,\ell'\in\{1,\dots ,p-1\}$,
\[
  \bigl(M_1^{\ell},M_2^{\rho^{r_2}(\ell)},\dots ,M_k^{\rho^{r_k}(\ell)}\bigr)
  =
  \bigl(M_1^{\ell'},M_2^{\rho^{r_2'}(\ell')},\dots ,M_k^{\rho^{r_k'}(\ell')}\bigr).
\]
The construction never rotates the first block, and equality holds
componentwise, so $M_1^{\ell}=M_1^{\ell'}$ and hence $\ell=\ell'$. Comparing the
$i$-th components gives $\rho^{r_i}(\ell)=\rho^{r_i'}(\ell)$, so
$\rho^{\,r_i-r_i'}$ fixes $\ell$. Since $\rho$ is a $(p-1)$-cycle, every non-identity power of a $(p-1)$ cycle is fixed-point-free. Therefore, $r_i\equiv r_i'\pmod{p-1}$; as $r_i,r_i'\in\{0,\dots ,p-2\}$ this
gives $r_i=r_i'$, a contradiction.
\end{proof}

\begin{lemma}\label{lem:recover}
Let $M=(M_1,\dots ,M_k)$ be a perfect matching of $L_{\V}$ in which
$M_2,\dots ,M_k$ are $p$-cycles. Suppose the cycles of $M$ are divided into
$k$ groups so that the cycle lengths in each group sum to $p$. Then the groups
are exactly the cycles of $M$ lying in $V_1,\dots ,V_k$. In particular, $M$
determines $\V$.
\end{lemma}

\begin{proof}
A $p$-cycle already has length $p$, so a group containing one contains nothing
else. Hence each of the $k-1$ cycles of $M_2,\dots ,M_k$ is a group on its own.
The only cycles left are those of $M_1$, and they must make up the last group.
Each of these groups is the set of cycles inside one block $V_i$.
\end{proof}

\begin{lemma}\label{lem:diff-config}
Let $\V$ and $\V'$ be vertex configurations with $\V\ne\V'$ as partitions of
$[n]$. From Construction \ref{construction}, create 2 perfect sets $R^{\V}_{\mathbf r}(P_1,\dots ,P_k)$ and
$R^{\V'}_{\mathbf r'}(P_1',\dots ,P_k')$ from $\V$ and $\V'$. Then
$R^{\V}_{\mathbf r}(P_1,\dots ,P_k)$ and
$R^{\V'}_{\mathbf r'}(P_1',\dots ,P_k')$ are disjoint perfect sets.
\end{lemma}

\begin{proof}
Suppose $M$ lies in both perfect sets. The blocks of $\V'$ group the cycles of $M$ into $k$ groups, each with total length $p$. $M$ is a matching from the $\V$-construction, so its components in $V_2,...,V_k$ are $p$-cycles. By Lemma~\ref{lem:recover}, this grouping is the block partition of $\V$. Therefore, $\V=\V'$.
\end{proof}

\subsection{Counting $\mathcal{L}_{p,1;k}$ perfect sets}\label{subsec:count-configs}

\begin{lemma}\label{lem:count-sets}
For a fixed vertex partition $\V$, Construction~\ref{construction} produces
$|\mathcal P(L_{p,1})|((p-2)!)^{k-1}(p-1)^{k-1}$ pairwise disjoint perfect sets of $L_{\V}$.
\end{lemma}

\begin{proof}
There are $|\mathcal P(L_{p,1})|((p-2)!)^{k-1}$ choices of
$(P_1,\dots ,P_k)$ and $(p-1)^{k-1}$ choices of $\mathbf r$. Each choice gives a
perfect set by Lemma~\ref{lem:rotation-is-perfect}. Two choices with different
tuples give disjoint sets by Lemma~\ref{lem:diff-P}, and two with the same
tuple and different $\mathbf r$ give disjoint sets by Lemma~\ref{lem:diff-r}.
\end{proof}

\begin{lemma}\label{lem:count-configs}
The number of partitions of $[n]$ into $k$ unordered blocks of $p$ vertices each
is
\[
  \frac{n!}{(p!)^{k}\,k!}.
\]
\end{lemma}

\begin{proof}
This counts the ways to place $n$ distinguishable objects into $k$
indistinguishable boxes of $p$ elements each. Treating the boxes as
distinguishable first,
\begin{align*}
  \binom{pk}{p}\binom{pk-p}{p}\cdots\binom{p}{p}
  &= \frac{(pk)!}{p!\,(pk-p)!}\cdot\frac{(pk-p)!}{p!\,(pk-2p)!}\cdots
     \frac{p!}{p!\,0!}= \frac{(pk)!}{(p!)^{k}} .
\end{align*}
Dividing by $k!$ to correct for the ordering of the boxes gives
$\dfrac{(pk)!}{(p!)^{k}k!}=\dfrac{n!}{(p!)^{k}k!}$.
\end{proof}

\section{From the block graphs to \texorpdfstring{$L_{n,1}$}{L(n,1)}}
\label{sec:layers}

The graphs $L_{\V}$ of Section~\ref{sec:blocks} omit every edge that crosses
between blocks. To reach $L_{n,1}$ we must also handle those crossing edges, and
to keep the resulting perfect sets disjoint we must control which matchings they
use.

Throughout, $p$ is prime, $k\ge2$, $n=pk$, and $\sigma=(1\,2\,\cdots\,k)$ is the
cyclic permutation of the block indices $\{1,\dots ,k\}$. We write $V(x)$ for the
block containing $x$. A vertex configuration $\V=(V_1,\dots ,V_k)$ carries a
fixed indexing of its blocks, and this indexing is what $\sigma$ acts on.

\subsection{Layer graphs}\label{subsec:layers}

\begin{defn}\label{def:layer}
Let $\V$ be a vertex configuration and fix $\kappa\in\{1,\dots ,k-1\}$. The
\emph{layer graph} $Q^{\kappa}_{\V}$ is the spanning subgraph of $L_{n,1}$ whose
edge set is a complete bipartite graph $K_{p,p}$ on each pair
$V_i,V_{\sigma^{\kappa}(i)}'$. Since $\kappa\not\equiv0\pmod k$ we have
$\sigma^{\kappa}(i)\neq i$, so no edge $xx'$ occurs and $Q^{\kappa}_{\V}$ is
indeed a subgraph of $L_{n,1}$. It is $p$-regular with $kp^{2}=pn$ edges.
\end{defn}

\begin{example}\label{ex:layers} 
Continuing Example~\ref{ex:LV} and Example~\ref{ex:rotation}, take $L_{\V}\in\Lcal_{3,1;3}$ with blocks
$V_1=\{1,5,6\}$, $V_2=\{2,8,9\}$ and $V_3=\{3,4,7\}$. The graph $Q^1_{\V}$ places
a $K_{3,3}$ on each pair $V_i/V_{\sigma(i)}'$:
\[
Q^1_\V \;=\;
\begin{tikzpicture}[baseline=(current bounding box.center)]
  \tworows{1,5,6,2,8,9,3,4,7}{2',8',9',3',4',7',1',5',6'}
  \draw[ed] (t0)--(b0) (t0)--(b1) (t0)--(b2) (t1)--(b0) (t1)--(b1)
            (t1)--(b2) (t2)--(b0) (t2)--(b1) (t2)--(b2);
  \draw[ed] (t3)--(b3) (t3)--(b4) (t3)--(b5) (t4)--(b3) (t4)--(b4)
            (t4)--(b5) (t5)--(b3) (t5)--(b4) (t5)--(b5);
  \draw[ed] (t6)--(b6) (t6)--(b7) (t6)--(b8) (t7)--(b6) (t7)--(b7)
            (t7)--(b8) (t8)--(b6) (t8)--(b7) (t8)--(b8);
  \blocklabels{$V_1/V_2'$}{$V_2/V_3'$}{$V_3/V_1'$}
\end{tikzpicture}
\]
and $Q^2_{\V}$ places a $K_{3,3}$ on each pair $V_i/V_{\sigma^2(i)}'$:
\[
Q^2_\V \;=\;
\begin{tikzpicture}[baseline=(current bounding box.center)]
  \tworows{1,5,6,2,8,9,3,4,7}{3',4',7',1',5',6',2',8',9'}
  \draw[ed] (t0)--(b0) (t0)--(b1) (t0)--(b2) (t1)--(b0) (t1)--(b1)
            (t1)--(b2) (t2)--(b0) (t2)--(b1) (t2)--(b2);
  \draw[ed] (t3)--(b3) (t3)--(b4) (t3)--(b5) (t4)--(b3) (t4)--(b4)
            (t4)--(b5) (t5)--(b3) (t5)--(b4) (t5)--(b5);
  \draw[ed] (t6)--(b6) (t6)--(b7) (t6)--(b8) (t7)--(b6) (t7)--(b7)
            (t7)--(b8) (t8)--(b6) (t8)--(b7) (t8)--(b8);
  \blocklabels{$V_1/V_3'$}{$V_2/V_1'$}{$V_3/V_2'$}
\end{tikzpicture}
\]
Note that $L_{\V}\cup Q^1_{\V}\cup Q^2_{\V}=L_{9,1}$.
\end{example}

By repeated application of Hall's marriage theorem, $Q^{\kappa}_{\V}$ decomposes
into $p$ edge-disjoint perfect matchings. Since every edge of the layer joins
$V_i$ to $V'_{\sigma^{\kappa}(i)}$ with $\sigma^{\kappa}(i)\ne i$, every perfect
matching of $Q^{\kappa}_{\V}$ is a derangement of $[n]$.

We extend the terminology of Section~\ref{sec:blocks} to these spanning
subgraphs: a \emph{perfect set of $Q^{\kappa}_{\V}$} is a $1$-factorization of
$Q^{\kappa}_{\V}$, necessarily into $p$ perfect matchings, and a perfect set of
$L_{\V}$ is a $1$-factorization of $L_{\V}$ into $p-1$ perfect matchings. Since
$(p-1)+(k-1)p=n-1$, a perfect set of $L_{\V}$ together with one perfect set of
each layer $Q^{1}_{\V},\dots ,Q^{k-1}_{\V}$ (all pairwise disjoint) has exactly $n-1$ members; by
Lemma~\ref{lem:L-vs-Q} below it is therefore a perfect set of $L_{n,1}$.

\begin{defn}\label{def:matching-notation}
Fix $\kappa$ and $\V$, and let $M$ be a perfect matching of $Q^{\kappa}_{\V}$.
\begin{enumerate}
  \item For $x \in V_i$ and $y' \in V_{\sigma^\kappa(i)}'$ we write $M(x)=y$ if
        and only if $xy'\in M$, so that $M$ is a permutation of $[n]$.
  \item We apply $M$ notation to blocks as well, writing $M(V_i)=V_{\sigma^{\kappa}(i)}$,
        and more generally
        \begin{equation}\label{eq:shift}
          M^{t}(V_i)=V_{\sigma^{\kappa t}(i)}\qquad(t\ge0).
        \end{equation}
  \item $M^{s}$ denotes the $s$-fold composite of the permutation $M$ with
        itself, so $M^{s}(x)$ is the vertex reached from $x$ by following the
        matching $s$ times; thus $M^{s}\circ M^{s'}=M^{s+s'}$ for all
        $s,s'\ge0$. We write $M^{s}|_{W}$ for the restriction of $M^{s}$ to
        $W\subseteq[n]$; when $M^{s}(W)=W$ this restriction is a permutation of
        $W$, and then $(M^{s}|_{W})^{r}=M^{sr}|_{W}$ for all $r\ge0$.
\end{enumerate}
\end{defn}

Powers of $\sigma$ split it into disjoint cycles of equal length: $\sigma^\kappa$
has $\gcd(\kappa,k)$ cycles, each of length $k/\gcd(k,\kappa)$. These two
parameters are used frequently below, so we give them the canonical names, $d$ and $m$:
\begin{equation}\label{eq:dm}
  d=\gcd(\kappa,k),\qquad m=k/d .
\end{equation}

\begin{defn}\label{def:m-orbit}
Fix $p,k$ and a layer $Q_\V^\kappa$, and let $d,m$ be as in \eqref{eq:dm}. Write
\[
  \sigma^{\kappa}=(i_{1,1}\ \cdots\ i_{1,m})(i_{2,1}\ \cdots\ i_{2,m})\cdots
                  (i_{d,1}\ \cdots\ i_{d,m}).
\]
An \emph{$m$-orbit} is the index set $O_j=\{i_{j,1},\dots ,i_{j,m}\}$ of one of
these $d$ cycles. We use the same name for the corresponding collection of
blocks $\{V_i:i\in O_j\}$. The $d$ $m$-orbits partition $\{1,\dots ,k\}$ and each has exactly $m$ elements.
Since $\kappa\not\equiv0\pmod k$ we have $d<k$, hence $m\ge2$.
\end{defn}

\subsection{Tautness}\label{subsec:taut}

\begin{lemma}\label{lem:structure}
Let $M$ be a perfect matching of $Q^{\kappa}_{\V}$ and $x\in[n]$.
\begin{enumerate}
\item[(a)] $M^{t}(x)\in V(x)$ if and only if $m\mid t$. In particular $M^{m}$
      maps each block to itself: $M^m(V_i)=V_i$ for all $i$.
\item[(b)] The cycle of $M$ through $x$ has length $m\ell$, where $\ell$ is the
      length of the cycle of $M^{m}$ through $x$. Moreover $\ell\le p$, so no
      cycle of $M$ is longer than $mp=n/d$.
\item[(c)] The following are equivalent:
      \emph{(i)} every cycle of $M^{m}$ has length $p$;
      \emph{(ii)} every cycle of $M$ has length $n/d$;
      \emph{(iii)} $M$ has exactly $d$ cycles.
\end{enumerate}
\end{lemma}

\begin{proof}
(a) Let $x\in V_i$. By \eqref{eq:shift}, $M^{t}(x)\in V_i$ if and only if
$\sigma^{\kappa t}(i)=i$. Since $\sigma$ is a $k$-cycle, $\sigma^{s}$ fixes an
index if and only if $k\mid s$, so the condition is $k\mid\kappa t$. Write
$\kappa=d\kappa_0$; since $\gcd(k,\kappa)=d$ we have $\gcd(\kappa_0,m)=1$, so
the condition becomes $m\mid\kappa_0t$ and hence $m\mid t$. Taking $t=m$ gives
$M^{m}(V_i)=V_i$.

(b) Let $L\ge1$ be least with $M^{L}(x)=x$, and $\ell\ge1$ least with
$(M^{m})^{\ell}(x)=x$. Since $M^{L}(x)=x\in V(x)$, part (a) gives $m\mid L$, say
$L=ms$. Then $(M^{m})^{s}(x)=x$, so $\ell\le s$ by minimality of $\ell$.
Conversely $M^{m\ell}(x)=(M^{m})^{\ell}(x)=x$, so $L\le m\ell$ by minimality of
$L$; as $L=ms$ this gives $s\le\ell$. Hence $s=\ell$ and $L=m\ell$. Finally, by
(a) the cycle of $M^{m}$ through $x$ lies in the $p$-element set $V(x)$, so
$\ell\le p$ and $L\le mp=n/d$.

(c) (i)$\Rightarrow$(ii) is immediate from (b). For (ii)$\Rightarrow$(iii), the
cycles of $M$ partition $[n]$, so there are $n/(n/d)=d$ of them. For
(iii)$\Rightarrow$(ii), $d$ cycles partition $[n]$ and each has at most $n/d$
vertices by (b), so each has exactly $n/d$. For (ii)$\Rightarrow$(i),
$m\ell=mp$ forces $\ell=p$ for every $x$.
\end{proof}

\begin{defn}\label{def:taut}
Call a perfect matching of $Q^{\kappa}_{\V}$ \emph{taut} if it satisfies the
equivalent conditions of Lemma~\ref{lem:structure}(c). For $d=1$, tautness says
exactly that $M$ is an $n$-cycle, where $n=pk$. For $d \geq 2$, a taut matching $M$ consists of $d$ disjoint cycles, each of length $mp$.
\end{defn}

\begin{remark}[Powers and return maps]\label{rem:powers}
Three consequences of Lemma~\ref{lem:structure} are used repeatedly.
\begin{enumerate}
\item[(i)] $M^{m}$ maps each block $V_i$ onto itself, so it restricts to a
      permutation $M^{m}|_{V_i}$ of the $p$-element set $V_i$; every cycle of
      $M^{m}$ lies inside a single block. We call $M^{m}|_{V_i}$ the
      \emph{return map at $V_i$}.
\item[(ii)] $M$ is taut if and only if $M^{m}|_{V_i}$ is a $p$-cycle for every
      $i$; equivalently, the cycle of $M^{m}$ through any $x$ is all of $V(x)$.
\item[(iii)] If $M$ is taut then every cycle of $M^{m}$ has length $p$, so
      $M^{mp}$ is the identity on $[n]$ and $M^{ms}$ depends only on
      $s\bmod p$. Exponents of the form $M^{ms}$ are therefore read modulo $p$
      whenever $M$ is taut.
\end{enumerate}
\end{remark}

\begin{lemma}\label{lem:rigid}
Fix a prime $p$ and let $k=n/p$. Suppose $M$ is taut as a matching of
$Q^{\kappa}_{\V}$ relative to a vertex configuration $\V$, and also taut as a
matching of $Q^{\kappa'}_{\V'}$ relative to a vertex configuration $\V'$. Then
$\gcd(\kappa,k)=\gcd(\kappa',k)$, $\V=\V'$ as partitions of $[n]$, and
$Q^{\kappa}_{\V}=Q^{\kappa'}_{\V'}$ as graphs. If in addition the common blocks
are indexed the same way in $\V$ and $\V'$, then $\kappa=\kappa'$.
\end{lemma}

\begin{proof}
Let $d,m$ and $d',m'$ be the quantities \eqref{eq:dm} for the two layers. By
Lemma~\ref{lem:structure}(c)(ii), every cycle of $M$ has length $n/d$ and also
length $n/d'$, so $d=d'$ and $m=m'$.

Let $V$ be a block of $\V$. By Lemma~\ref{lem:structure}(a), $M^{m}$ maps $V$
onto itself, so $V$ is a union of cycles of $M^{m}$; each such cycle has length
$p=|V|$ by (c)(i), so $V$ is a single cycle of $M^{m}$. As the blocks of $\V$
and the cycles of $M^{m}$ both partition $[n]$, the blocks of $\V$ are exactly
the cycles of $M^{m}$. The same argument applies to $\V'$ with the same exponent
$m$, so $\V=\V'$. Both layers therefore have the same blocks, and by
\eqref{eq:shift} the edge set of each layer is determined by the map
$V\mapsto M(V)$ on those blocks, which is the same map in both cases; hence
$Q^{\kappa}_{\V}=Q^{\kappa'}_{\V'}$.

Finally, suppose the blocks are indexed the same way in $\V$ and $\V'$. Then
\eqref{eq:shift} in each layer gives
$V_{\sigma^{\kappa}(i)}=M(V_i)=V_{\sigma^{\kappa'}(i)}$ for all $i$, so
$\sigma^{\kappa}=\sigma^{\kappa'}$ and $\kappa\equiv\kappa'\pmod k$. Since
$\kappa,\kappa'\in\{1,\dots ,k-1\}$, this gives $\kappa=\kappa'$.
\end{proof}

\begin{remark}\label{rem:rigid-caveat}
The final clause cannot be dropped. A taut $M$ determines the partition and the
induced block permutation $V\mapsto M(V)$, but not the labelling of the blocks:
relabelling the blocks so that the same block permutation reads as a different
power of $\sigma$ changes $\kappa$ while leaving the graph $Q^{\kappa}_{\V}$
unchanged. What is intrinsic to $M$, and what is used below, is the partition
and the layer graph itself.
\end{remark}

\begin{cor}\label{cor:disjoint}
For a fixed block size $p$, a taut matching determines both the partition $\V$
and the layer graph containing it. Consequently, taut matchings arising from
distinct vertex configurations are distinct, and for a fixed configuration taut
matchings of $Q^{\kappa}_{\V}$ and $Q^{\kappa'}_{\V}$ with $\kappa\ne\kappa'$
are distinct.
\end{cor}

\begin{proof}
Immediate from Lemma~\ref{lem:rigid}.
\end{proof}

\begin{remark}\label{rem:fixed-p}
The prime $p$ must be fixed here. The blocks are recovered as the cycles of
$M^{k/d}$ with $k=n/p$, so a different block size uses a different power of $M$
and may well recover a different partition.
\end{remark}

Blocks in a common $\sigma^{\kappa}$-orbit are related by $M$, and the next
lemma says that following the matching one step carries the cycles of $M^{m}$ in
one block onto those in the next, without distortion.

\begin{lemma}\label{lem:transport}
Let $M$ be a perfect matching of $Q^{\kappa}_{\V}$.
\begin{enumerate}
\item[(a)] For all $s,s'\ge0$ and all $x\in[n]$,
      $M^{s}\bigl(M^{s'}(x)\bigr)=M^{s'}\bigl(M^{s}(x)\bigr)=M^{s+s'}(x)$. In
      particular, if the return map at $V_i$ is
      $M^{m}|_{V_i}=(x_1\ x_2\ \cdots\ x_p)$, then
      $M^{ms}|_{V_i}=\bigl(M^{m}|_{V_i}\bigr)^{s}=(x_1\ x_2\ \cdots\ x_p)^{s}$.
\item[(b)] For every $x$, the permutation $M$ maps the cycle of $M^{m}$ through
      $x$ bijectively onto the cycle of $M^{m}$ through $M(x)$; in particular the
      two cycles have the same length.
\item[(c)] For every $i$, the return map at $V_i$ is a $p$-cycle if and only if
      the return map at $V_{\sigma^{\kappa}(i)}$ is. Hence if the return map at
      one block is a $p$-cycle, the same holds at every block of that block's
      $\sigma^{\kappa}$-orbit.
\end{enumerate}
\end{lemma}

\begin{proof}
(a) The first claim is immediate from
$M^{s}\circ M^{s'}=M^{s+s'}=M^{s'}\circ M^{s}$. For the second,
Lemma~\ref{lem:structure}(a) gives $M^{m}(V_i)=V_i$, so
restricting to $V_i$ commutes with taking powers:
$M^{ms}|_{V_i}=\bigl(M^m|_{V_i}\bigr)^{s}$, which in cycle notation is
$(x_1\ x_2\ \cdots\ x_p)^{s}$.

(b) The cycle of $M^{m}$ through $x$ is $\{M^{ms}(x):s\ge0\}$ and the cycle
through $M(x)$ is $\{M^{ms}(M(x)):s\ge0\}$. By (a),
$M\bigl(M^{ms}(x)\bigr)=M^{ms}\bigl(M(x)\bigr)$ for all $s$, so $M$ restricts to
a bijection between the two sets; in particular they have the same size.

(c) By Remark~\ref{rem:powers}(i), each cycle of $M^{m}$ lies inside a single
block, and by \eqref{eq:shift} we have $M(V_i)=V_{\sigma^{\kappa}(i)}$. So (b)
matches the cycles of $M^{m}$ lying in $V_i$ with those lying in
$V_{\sigma^{\kappa}(i)}$, one to one and preserving lengths; in particular the
two blocks contain the same number of cycles of $M^{m}$. Each block has $p$
elements, so its return map is a $p$-cycle exactly when it contains a single
cycle of $M^{m}$, and these counts agree. Iterating along the orbit gives the
last statement.
\end{proof}

\begin{example}\label{ex:m-orbit}
Let $L_{\V}\in\Lcal_{3,1;6}$ with $V_1=\{v_1,v_3,v_5\}$, and consider the layer
$Q^2_{\V}$. Here $(1\ 2\ 3\ 4\ 5\ 6)^2=(1\ 3\ 5)(2\ 4\ 6)$, so $d=2$ and $m=3$,
and the two $m$-orbits are $\{1,3,5\}$ and $\{2,4,6\}$. A matching $M$ of
$Q^2_{\V}$ therefore shifts blocks as
\[
  M(V_1)=V_3,\quad M(V_3)=V_5,\quad M(V_5)=V_1;\qquad
  M(V_2)=V_4,\quad M(V_4)=V_6,\quad M(V_6)=V_2 .
\]
Since $M^3(V_1)=V_1$, the return map $M^{3}|_{V_1}$ is a permutation of
$\{v_1,v_3,v_5\}$; say it is the $3$-cycle $(v_1\ v_3\ v_5)$, so $M^3(v_1)=v_3$.
Then the cycle of $M$ through $v_1$ has length $m\ell=3\cdot3=9=n/d$, and $M$ is
taut (Definition~\ref{def:taut}) precisely when the return maps at $V_1$ and at
$V_2$ are both $3$-cycles, by Lemma~\ref{lem:transport}(c).
\end{example}

\subsection{Twisting}\label{subsec:twisting}

Applying the rotation vector $\mathbf r$ of Section~\ref{sec:blocks} to a layer
produces non-taut matchings, so we introduce a related but different device.

\begin{defn}\label{def:twistvector}
Consider $Q^{\kappa}_{\V}$ with $d=\gcd(k,\kappa)$ and $m=k/d$, and write
\[
  \sigma^\kappa=(i_{1,1}\ \cdots\ i_{1,m})(i_{2,1}\ \cdots\ i_{2,m})\cdots
                (i_{d,1}\ \cdots\ i_{d,m}).
\]
A \emph{twist vector} is a vector $\mathbf a=(a_1,\dots ,a_k)$ with
$a_i \in \Z_p^{\times}$ such that
\begin{equation}\label{eq:A}
  \sum_{i\in O}a_i\equiv0\pmod p
  \qquad\text{for every $m$-orbit $O$.}
\end{equation}
\end{defn}

\begin{remark}\label{rem:per-orbit}
Condition \eqref{eq:A} is imposed on each $m$-orbit separately: the sum rule
must hold within each orbit, not merely over all $k$ coordinates at once. See
Example~\ref{ex:twistvector}.
\end{remark}

\begin{example}\label{ex:twistvector}
Let $L_{\V}\in\Lcal_{3,1;4}$ with $\V=(V_1,V_2,V_3,V_4)$, and consider the layer
$Q_\V^1$. Since $\gcd(4,1)=1$ there is a single $m$-orbit, $(1\ 2\ 3\ 4)$, and
\[
  \mathbf a=(-1,1,-1,1)
\]
is admissible because $a_1+a_2+a_3+a_4\equiv0\pmod 3$.

Now consider $Q_\V^2$. Here $\gcd(4,2)=2$, so there are $d=4/2=2$ $m$-orbits, each with $m=4/2=2$
 blocks, namely $(1\ 2\ 3\ 4)^2=(1\ 3)(2\ 4)$. An admissible twist
vector is
\[
  \mathbf a=(-1,-1,1,1),
\]
since $a_1+a_3\equiv0$ and $a_2+a_4\equiv0 \pmod 3$. Note that $(-1,1,-1,1)$ is
\emph{not} admissible here, since $a_1+a_3\not\equiv0\pmod 3$ (violating the
condition on a single orbit already suffices).
\end{example}

\begin{lemma}\label{lem:exists-a}
A twist vector for $Q^{\kappa}_{\V}$ exists if and only if $p$ is odd, or $p=2$
and $m$ is even.
\end{lemma}

\begin{proof}
The $d$ orbits are disjoint and all of size $m$, so it suffices to produce a
list of $m$ units of $\Z_p$ summing to $0$. If $p=2$ the only unit is $1$, so
the sum is $m \bmod 2$ and $m$ must be even; conversely, if $p=2$ and $m$ is
even, then $(1,1,\dots ,1)$ works.

Now let $p\ge3$. If $m$ is even, take the sub-vector on each $m$-orbit to be
$(1,-1,1,-1,\dots)$. If $m$ is odd then $m\ge3$, since $m\ge2$ and the case
$m=2$ is even; take $(1,1,-2)$ followed by $(m-3)/2$ pairs $(1,-1)$. In each
case all entries are units of $\Z_p$ and they sum to $0$.
\end{proof}

\begin{defn}\label{def:twist}
Fix a twist vector $\mathbf a$ and let $M$ be a taut matching of
$Q^{\kappa}_{\V}$. For $t\in\Z_p$ define $M_t$ by
\begin{equation}\label{eq:twistdef}
  M_t(x)\;=\;M\bigl(M^{\,m a_i t}(x)\bigr)\qquad\text{for }x\in V_i,
\end{equation}
where the exponent $a_it$ is read modulo $p$. That is, before following the
matching out of $V_i$, first move $a_it$ steps along the return map at $V_i$.
The \emph{twist family} of $M$ is $\Fcal(M):=\{M_t:t\in\Z_p\}$.
\end{defn}

\begin{lemma}\label{lem:twist}
Let $\mathbf a$ be a fixed twist vector and $M$ a taut perfect matching of the
layer $Q^\kappa_{\V}$. Then:
\begin{enumerate}
\item[(a)] Each $M_t$ is a perfect matching of $Q^{\kappa}_{\V}$, and $M_0=M$.
\item[(b)] $M_t$ is taut, and moreover $M_t^{\,m}=M^{m}$.
\item[(c)] The $p$ matchings $M_t$ in $\Fcal(M)$ are pairwise edge-disjoint (so
      $|\Fcal(M)|=p$), and $(M_t)_{t'}=M_{t+t'}$.
\item[(d)] $\Fcal(M)$ is a perfect set of $Q^{\kappa}_{\V}$.
\item[(e)] Every taut matching of $Q^{\kappa}_{\V}$ lies in exactly one twist
      family formed with respect to $\mathbf a$; that is, the twist families
      partition the taut matchings of $Q^{\kappa}_{\V}$ into classes of size $p$.
\end{enumerate}
\end{lemma}

\begin{proof}
\textbf{(a)} Taking $t=0$ gives $M_0(x)=M(M^0(x))=M(x)$ for all $x$, so
$M_0=M$. Now fix $i$. Since $M^{\,ma_it}$ is a power of $M^m$ and $M^m$ maps
$V_i$ to itself, $M^{\,ma_it}|_{V_i}$ is a permutation of $V_i$. Composing it
with the bijection $M|_{V_i}:V_i\to V_{\sigma^{\kappa}(i)}$ shows that $M_t$
maps $V_i$ bijectively onto $V_{\sigma^{\kappa}(i)}$ for every $i$. Since the
blocks partition $[n]$ and $\sigma^{\kappa}$ permutes their indices, the images
$V_{\sigma^{\kappa}(i)}$ also partition $[n]$, so $M_t$ is a permutation of
$[n]$. Finally, each edge $x\,M_t(x)'$ joins $V_i$ to $V'_{\sigma^{\kappa}(i)}$
and is therefore an edge of $Q^{\kappa}_{\V}$.

\medskip
\textbf{(b)} Fix an $m$-orbit and denote its blocks by $V_{i_0},\dots ,V_{i_{m-1}}$
with $i_{r+1}=\sigma^{\kappa}(i_r)$. Since $M$ is taut, the return map
$M^{m}|_{V_{i_0}}$ is a $p$-cycle, which we write as
$(x_1\ x_2\ \cdots\ x_p)$. Let $x\in V_{i_0}$ be arbitrary. We show by induction
that for $0\le r\le m$,
\begin{equation}\label{eq:taut-induction}
  M_t^{\,r}(x)=M^{\,r+mt(a_{i_0}+\cdots+a_{i_{r-1}})}(x).
\end{equation}
For $r=0$ both sides are $x$: the left side is $M_t^{0}(x)=x$, and on the right
the sum $a_{i_0}+\cdots+a_{i_{r-1}}$ is empty, hence $0$, so the exponent is $0$.

Suppose \eqref{eq:taut-induction} holds for some $r<m$ and set
$x_r:=M_t^{\,r}(x)\in V_{i_r}$. By definition of the twist,
\[
  M_t(x_r)=M\bigl(M^{\,ma_{i_r}t}(x_r)\bigr)=M^{\,ma_{i_r}t+1}(x_r),
\]
so composing with \eqref{eq:taut-induction} and using
Lemma~\ref{lem:transport}(a) to merge the two powers of $M$ applied to $x$
(both are powers of $M$ landing at the same intermediate point, so the exponents
add),
\[
  M_t^{\,r+1}(x)=M^{\,ma_{i_r}t+1}\bigl(M^{\,r+mt(a_{i_0}+\cdots+a_{i_{r-1}})}(x)\bigr)
  =M^{\,(r+1)+mt(a_{i_0}+\cdots+a_{i_r})}(x),
\]
which is \eqref{eq:taut-induction} at $r+1$.

Taking $r=m$ closes the $m$-orbit back to $V_{i_0}$ and gives, again by
Lemma~\ref{lem:transport}(a) since
$M^{m(1+tS_O)}|_{V_{i_0}}=\bigl(M^m|_{V_{i_0}}\bigr)^{1+tS_O}$,
\[
  M_t^{\,m}(x)=M^{\,m(1+tS_O)}(x)=(x_1\ x_2\ \cdots\ x_p)^{\,1+tS_O}(x),
  \qquad S_O:=a_{i_0}+\cdots+a_{i_{m-1}} .
\]
Since $\mathbf a$ is a twist vector, $S_O\equiv0\pmod p$, so
$1+tS_O\equiv1\pmod p$ for every $t$, and by Remark~\ref{rem:powers}(iii) the
exponent may be read modulo $p$. Hence
$M_t^{\,m}|_{V_{i_0}}=(x_1\ x_2\ \cdots\ x_p)=M^{m}|_{V_{i_0}}$, a $p$-cycle,
for every $t$.

The block $V_{i_0}$ was an arbitrary member of an arbitrary $m$-orbit and every
block lies in some $m$-orbit, so $M_t^{\,m}=M^{m}$ on all of $[n]$ and every
return map of $M_t$ is a $p$-cycle; thus $M_t$ is taut by
Remark~\ref{rem:powers}(ii). (Alternatively, one may verify the $p$-cycle
condition at a single block of each $m$-orbit and propagate it along the orbit
using Lemma~\ref{lem:transport}(c), since $M_t$ is a matching of
$Q^{\kappa}_{\V}$ by (a).)

\medskip
\textbf{(c)} Let $t\ne t'$ in $\Z_p$ and let $x\in V_i$. If $M_t(x)=M_{t'}(x)$
then, since $M$ is injective,
\[
  M^{ma_it}(x)=M^{ma_it'}(x),\qquad\text{that is,}\qquad
  \bigl(M^{m}|_{V_i}\bigr)^{a_i(t-t')}(x)=x .
\]
Now $a_i\in\Z_p^{\times}$ and $t-t'\not\equiv0\pmod p$, so
$a_i(t-t')\not\equiv0\pmod p$. Since $p$ is prime, every non-identity power of
the $p$-cycle $M^{m}|_{V_i}$ is again a $p$-cycle and so has no fixed point in
$V_i$, contradicting the displayed equation. Hence $M_t(x)\ne M_{t'}(x)$ for
\emph{every} vertex $x$, so $M_t$ and $M_{t'}$ share no edge; in particular the
$p$ matchings $M_t$ are distinct and $|\Fcal(M)|=p$.

By (b), $M_t$ is taut with $M_t^{\,m}=M^{m}$, so twisting it is governed by the
same powers: for $x\in V_i$,
\[
  (M_t)_{t'}(x)=M_t\bigl(M_t^{\,ma_it'}(x)\bigr)=M_t\bigl(M^{\,ma_it'}(x)\bigr)
  =M\bigl(M^{\,ma_it}\bigl(M^{\,ma_it'}(x)\bigr)\bigr)=M_{t+t'}(x),
\]
using that $M^{\,ma_it'}(x)$ again lies in $V_i$.

\medskip
\textbf{(d)} By (a) each $M_t$ is a perfect matching of $Q^{\kappa}_{\V}$, and by
(c) the $p$ of them are pairwise edge-disjoint, so we check that they cover
the layer. Fix $x\in V_i$. The return map $M^m|_{V_i}$ is a $p$-cycle
$(x_1\ x_2\ \cdots\ x_p)$, and by Lemma~\ref{lem:transport}(a),
$M^{ma_it}|_{V_i}=(x_1\ x_2\ \cdots\ x_p)^{a_it}$. Since $a_i$ is a unit and $p$
is prime, the exponents $a_it$ run over all of $\Z_p$ as $t$ does, and distinct
powers of a $p$-cycle move $x$ to distinct vertices; hence the $p$ vertices
$M^{ma_it}(x)$, $t\in\Z_p$, are distinct and exhaust $V_i$. Applying the
bijection $M|_{V_i}$ shows that $x$ is matched to each vertex of
$V_{\sigma^{\kappa}(i)}'$ exactly once as $t$ varies, so every layer edge at $x$
lies in exactly one $M_t$. Since every edge of the layer has exactly one endpoint
in the unprimed copy of $[n]$, every edge lies in exactly one member of
$\Fcal(M)$, and the layer has $pn$ edges and each $M_t$ has
$n$. Hence $\Fcal(M)$ is a perfect set of $Q^{\kappa}_{\V}$.

\medskip
\textbf{(e)} The twist vector $\mathbf a$ is fixed throughout. For $x\in V_i$,
the map $M^{\,ma_it}|_{V_i}=\bigl(M^{m}|_{V_i}\bigr)^{a_it}$ is a power of a
$p$-cycle by Lemma~\ref{lem:transport}(a), so it depends only on $t\bmod p$, and
hence so does $M_t$; we therefore index $\Fcal(M)=\{M_t\}$ by $t\in\Z_p$.

Let $N$ be any taut matching of $Q^{\kappa}_{\V}$. Since $N$ is taut, $\Fcal(N)$
is defined and $N=N_0\in\Fcal(N)$ by (a), so $N$ lies in at least one twist
family. For uniqueness, suppose $N\in\Fcal(M)$ for some taut $M$, say $N=M_s$.
By (b), $M_s$ is taut, so we may twist it, and by (c),
\[
  N_t=(M_s)_t=M_{s+t}\qquad\text{for every }t\in\Z_p .
\]
As $t$ ranges over $\Z_p$ so does $s+t$, whence
\[
  \Fcal(N)=\{N_t: t\in\Z_p\}=\{M_{s+t}: t\in\Z_p\}=\Fcal(M).
\]
So every twist family containing $N$ equals $\Fcal(N)$, and $N$ lies in exactly
one twist family. Each family has exactly $p$ members by (c).
\end{proof}

\subsection{Counting taut matchings and perfect sets of a layer}
\label{subsec:layer-count}

\begin{lemma}\label{lem:count-taut}
Each layer $Q^{\kappa}_{\V}$ has exactly $(p!)^{k}/p^{\,d}$ taut perfect
matchings.
\end{lemma}

\begin{proof}
Consider a single $m$-orbit with blocks $V_{i_0},V_{i_1},\dots ,V_{i_{m-1}}$,
where $i_{r+1}=\sigma^{\kappa}(i_r)$. On this orbit a perfect matching $M$ of the
layer restricts to a list of bijections
\[
  g_{i_r}:V_{i_r}\longrightarrow V_{i_{r+1}}\qquad
  (r=0,1,\dots ,m-1,\ \text{indices mod } m),
\]
where $g_{i_r}$ is the edge set placed on the pair $V_{i_r}/V_{i_{r+1}}'$; and
conversely any such list determines $M$ on the orbit.

Choose first the bijections $g_{i_1},\dots ,g_{i_{m-1}}$ on the pairs
$V_{i_1}/V_{i_2}'$, $V_{i_2}/V_{i_3}'$, \dots , $V_{i_{m-1}}/V_{i_0}'$; there
are $(p!)^{m-1}$ ways to do so. Leave $g_{i_0}$
undetermined for the moment, since it is the one that will be constrained.
Composing the chosen maps gives
\[
  \varphi:V_{i_1}\longrightarrow V_{i_0},\qquad
  \varphi=g_{i_{m-1}}\circ\cdots\circ g_{i_1},
\]
which is exactly $M^{\,m-1}|_{V_{i_1}}$; it is a bijection, being a composite of
bijections. By Remark~\ref{rem:powers}(ii) together with
Lemma~\ref{lem:transport}(c), $M$ is taut if and only if, for each $m$-orbit,
the return map at the single block $V_{i_0}$ of that orbit,
\[
  M^{m}|_{V_{i_0}}=g_{i_{m-1}}\circ\cdots\circ g_{i_1}\circ g_{i_0}
  =\varphi\circ g_{i_0},
\]
is a $p$-cycle. Here $\varphi$ is already fixed, so $g_{i_0}\mapsto\varphi\circ
g_{i_0}$ is a bijection from the $p!$ bijections $V_{i_0}\to V_{i_1}$ onto the
$p!$ permutations of $V_{i_0}$; equivalently $g_{i_0}=\varphi^{-1}\circ c$, where
$c$ is the desired return map. So the admissible choices of $g_{i_0}$ correspond
one for one to the $p$-cycles $c$ on $V_{i_0}$. (It is not the case that
$g_{i_0}$ is itself a $p$-cycle: for $m\ge2$ it is a bijection between two
different blocks, so the term does not apply to it.)

There are $(p-1)!$ $p$-cycles, so the number of admissible restrictions to one
$m$-orbit is
\[
  (p!)^{m-1}\,(p-1)!=(p!)^{m-1}\cdot\frac{p!}{p}=\frac{(p!)^{m}}{p} .
\]
The $d$ orbits partition the blocks, the restrictions of $M$ to different orbits
are chosen independently, and tautness is one condition per orbit, so the total
is
\[
  \frac{(p!)^{m d}}{p^d} = \frac{(p!)^{k}}{p^d}. \qedhere
\]
\end{proof}

\begin{thm}\label{thm:main}
Let $p$ be prime, $k\ge2$, $n=pk$ and $\kappa\not\equiv0\pmod k$, and let
$d=\gcd(\kappa,k)$ and $m=k/d$. If $p$ is odd or $m$ is even, then for each twist
vector $\mathbf a$ the twist families $\Fcal(M)$ partition the taut matchings of
$Q^{\kappa}_{\V}$ into
\[
  N_{\kappa}\;=\;\frac{(p!)^{k}}{p^{\,d+1}}
\]
disjoint perfect sets of $Q^{\kappa}_{\V}$, every member of which is taut.
\end{thm}

\begin{proof}
A twist vector $\mathbf a$ exists by Lemma~\ref{lem:exists-a}, and taut
matchings exist by Lemma~\ref{lem:count-taut}, since $(p!)^{k}/p^{\,d}\ge1$. By
Lemma~\ref{lem:twist}(e) the families $\Fcal(M)$ partition the set of taut
matchings into classes of size exactly $p$ by (c), and by (d) each class is a
perfect set of the layer. Hence the number of classes is
\[
  \frac{1}{p}\cdot\frac{(p!)^{k}}{p^{\,d}}=\frac{(p!)^{k}}{p^{\,d+1}} . \qedhere
\]
\end{proof}

\begin{cor}\label{cor:capacity}
Suppose $p$ is odd, or $p=2$ and $n$ is a power of $2$, and let $k=n/p\ge2$. Let
\[
  T=|\mathcal P(L_{p,1})|\,\bigl((p-2)!\bigr)^{k-1}(p-1)^{k-1}
\]
be the number of perfect sets of $L_{\V}$ from Lemma~\ref{lem:count-sets}. Then
$N_\kappa\ge T$ for every $\kappa\in\{1,\dots ,k-1\}$.
\end{cor}

\begin{proof}
Let $d=\gcd(\kappa,k)$ and $m=k/d\ge2$. If $p=2$ and $n=2^a$ then $m$ divides
$k=2^{a-1}$, so $m$ is even; hence in all cases Theorem~\ref{thm:main} applies and
$N_\kappa=(p!)^k/p^{\,d+1}$. Since $p!=p(p-1)(p-2)!$,
\[
  N_\kappa=p^{\,k-d-1}\bigl((p-1)!\bigr)^{k},\qquad
  T=|\mathcal P(L_{p,1})|\bigl((p-1)!\bigr)^{k-1}.
\]
As $d\le k/2$ we have $k-d-1\ge0$, so it suffices that
$|\mathcal P(L_{p,1})|\le(p-1)!$. This holds because
$|\mathcal P(L_{p,1})|\le d_p/(p-1)$ and $d_p\le p!/2\le(p-1)(p-1)!$.
\end{proof}

\subsection{Assembling the lower bound}\label{subsec:assembly}

\begin{lemma}\label{lem:L-vs-Q}
Fix a vertex configuration $\V$. The subgraphs
$L_{\V},Q^{1}_{\V},\dots ,Q^{k-1}_{\V}$ are pairwise edge-disjoint and their
edge union is $L_{n,1}$. Consequently no perfect matching of one of them is a
perfect matching of another; in particular, no perfect matching of $L_{\V}$ is a
perfect matching of any layer.
\end{lemma}

\begin{proof}
Let $xy'$ be an edge of $L_{n,1}$, so $x\ne y$; say $x\in V_i$ and $y\in V_j$. If
$i=j$ then the edge lies in $L_{\V}$, and it lies in no layer since every layer
edge joins $V_i$ to $V'_{\sigma^{\kappa}(i)}$ with $\sigma^{\kappa}(i)\ne i$. If
$i\ne j$ then, since $\sigma$ is a $k$-cycle, $j=\sigma^{\kappa}(i)$ for exactly
one $\kappa\in\{1,\dots ,k-1\}$, so the edge lies in exactly one layer. This
proves the decomposition, and a perfect matching of a subgraph uses only that
subgraph's edges, so it cannot be a perfect matching of an edge-disjoint one.

The two families are also separated by cycle type, which is what will be used
across different vertex configurations. Every cycle of a perfect matching of
$L_{\V}$ is contained in a single block, so has length at most $p$; whereas by
Lemma~\ref{lem:structure}(b) every cycle of a perfect matching of
$Q^{\kappa}_{\V}$ has length $m\ell$ with $m\ge2$, and for a taut matching has
length exactly $pm\ge2p$.
\end{proof}

\begin{remark}\label{rem:cycle-type-caveat}
The corresponding statement ``$M$ has $k$ cycles, each of length $p$'' holds when every block input is drawn from $\mathcal P(C_p)$. In Construction~\ref{construction} it is only guaranteed for blocks $2,...,k$, which is precisely what Lemma~\ref{lem:recover} requires.
\end{remark}

\begin{thm}\label{thm:construction}
Let $p$ be a prime dividing $n$ with $k=n/p\ge2$, and suppose $p$ is odd, or $p=2$
and $n$ is a power of $2$. Let
$T=|\mathcal P(L_{p,1})|\,((p-2)!)^{k-1}(p-1)^{k-1}$ be the number of perfect sets
of Lemma~\ref{lem:count-sets}. If $T\le N_\kappa$ for every
$\kappa\in\{1,\dots ,k-1\}$, then
\[
  \bigl|\calP(L_{n,1})\bigr|\;\ge\;\frac{n!}{(p!)^{k}k!}\,T.
\]
\end{thm}

\begin{proof}
Fix a partition $\V$ indexed as in Construction~\ref{construction}. Let
$A_1,\dots ,A_T$ be the perfect sets of $L_{\V}$ from Lemma~\ref{lem:count-sets},
and for each $\kappa$ choose pairwise disjoint perfect sets
$B^{\kappa}_1,\dots ,B^{\kappa}_T$ of $Q^{\kappa}_{\V}$, which is possible by
Theorem~\ref{thm:main} and $T\le N_\kappa$. Let
$S_j=A_j\cup B^1_j\cup\cdots\cup B^{k-1}_j$, and by Lemma~\ref{lem:L-vs-Q} each $S_j$
is a perfect set of $L_{n,1}$. For $j\neq j'$ we have $A_j\cap A_{j'}=\emptyset$ and
$B^{\kappa}_j\cap B^{\kappa}_{j'}=\emptyset$ for every $\kappa$ by construction, and a
matching of one of $L_\V,Q^1_\V,\dots,Q^{k-1}_\V$ is never a matching of another
by Lemma~\ref{lem:L-vs-Q}. Therefore $S_j\cap S_{j'}=\emptyset$.

Now let $\V\ne\V'$ with corresponding perfect sets and suppose a matching $M$ lies in a set built from each. By the cycle-length separation in Lemma~\ref{lem:L-vs-Q}, $M$ is a block matching in both configurations or a taut layer matching in both. If it
is a $Q$-layer matching, Lemma~\ref{lem:rigid} gives $\V=\V'$. If it is an $L$-part
matching, Lemma~\ref{lem:recover} shows the blocks of both $\V$ and $\V'$ are the
groups of cycles of $M$, so $\V=\V'$. This contradiction shows sets from different
partitions are disjoint.

Finally, counting partitions with
Lemma~\ref{lem:count-configs} gives the bound.
\end{proof}

We reproduce and prove the statement of Theorem \ref{thm:general_lower_bound} below.

\begin{restatethm}{Theorem \ref{thm:general_lower_bound}}
Let $n\ge 2$, let $p$ be the largest prime factor of $n$, set $k=n/p$, and let
$|\mathcal P(L_{p,1})|$ denote the maximum number of pairwise disjoint perfect sets of
$L_{p,1}$. Then
\[
  \lvert \calP(L_{n,1})\rvert
  \;\ge\;
  \frac{n!}{(p!)^{k}\, k!}\;|\mathcal P(L_{p,1})|\bigl((p-2)!\bigr)^{k-1} (p-1)^{\,k-1}.
\]
\end{restatethm}
\begin{proof}[Proof of Theorem~\ref{thm:general_lower_bound}]
Let $p$ be the largest prime factor of $n$ and $k=n/p$. If $k=1$ the bound reads
$|\calP(L_{n,1})|\ge|\mathcal P(L_{p,1})|$, which is the definition of
$\mathcal P(L_{p,1})$. Let $k\ge2$. If $p=2$ then $n$ has no odd prime factor, so
$n$ is a power of $2$; otherwise $p$ is odd. Then
Corollary~\ref{cor:capacity} gives $N_\kappa\ge T$ for every $\kappa$, so
Theorem~\ref{thm:construction} applies and gives
\[
  |\calP(L_{n,1})|\ge\frac{n!}{(p!)^kk!}\,T
  =\frac{n!}{(p!)^kk!}\,|\mathcal P(L_{p,1})|\bigl((p-2)!\bigr)^{k-1}(p-1)^{k-1}.\qedhere
\]
\end{proof}

Although our lower bound works for any prime satisfying our conditions, we state the bound only for the largest prime, since it yields the highest bound.

\begin{lemma}\label{prop:largest}
Let $p$ and $q$ be distinct primes dividing $n$, both admissible in
Theorem~\ref{thm:construction} (that is, odd, or $2$ with $n$ a power of $2$),
with $p>q$. Put $k=n/p$ and $k'=n/q$, and let $B_p,B_q$ denote the bounds found by Theorem~\ref{thm:construction} for $p$ and $q$ respectively. Then $B_p \geq B_q$.
\end{lemma}

\begin{proof}
If $n$ is a power of $2$ there is only one prime, so assume $n$ has an odd prime
factor. Then $2$ is not admissible, so $p>q\ge3$. Write $n=pqt$, so $k=n/p=qt$ and
$k'=n/q=pt$. Since $(p-2)!\,(p-1)=(p-1)!$ and $p!=p\,(p-1)!$, the bound with
the prime $p$ is
\[
  f(p)=\frac{n!}{(p!)^{k}k!}\,|\mathcal P(L_{p,1})|\,((p-1)!)^{k-1}=\frac{n!\,|\mathcal P(L_{p,1})|}{p^{k}\,(p-1)!\,k!}.
\]
Using $|\mathcal P(L_{p,1})|\ge(p-2)!$ gives $f(p)\ge\dfrac{n!}{(p-1)\,p^{k}\,k!}$. Using
$|\mathcal P(L_{q,1})|\le d_q/(q-1)\le q!/(2(q-1))$ gives
$f(q)\le\dfrac{n!}{2(q-1)\,q^{k'-1}\,k'!}$. So it suffices to show
\[
  \frac{(pt)!}{(qt)!}\;\ge\;\frac{(p-1)\,q}{2(q-1)}\Bigl(\frac{p^{q}}{q^{p}}\Bigr)^{t}.
\]
Since $x^{1/x}$ is decreasing for $x\ge e$, we have $q^{p}>p^{q}$, so the last
factor is at most $1$. As $q/(q-1)\le 3/2$, the right side is at most
$\tfrac34(p-1)<p$. The left side is a product of positive integers one of which
is $pt\ge p$, so it is at least $p$.
\end{proof}

\section{The cases \texorpdfstring{$n=5,6,7$}{n=5,6,7}}\label{sec:small}

This section proves Theorem~\ref{thm:small_cases}. Explicit analytical perfect partition
constructions for $L_{n,1}$ were difficult to find beyond the prime case and rotation construction, so we also modeled the search space amenable to computationally finding a perfect partition to push on smaller cases.

\subsection{The disjointness graph \texorpdfstring{$\Rcal_n$}{R(n)}}
\label{subsec:Rn}

\begin{defn}\label{def:Rn}
Let $R(n)$ denote the set of all reduced Latin squares of order $n$, equivalently
the set of all perfect sets of $L_{n,1}$. The graph $\Rcal_n$ has
$V(\Rcal_n)=R(n)$, with an edge between two perfect sets if and only if they
share no derangement:
\[
  E(\Rcal_n)=\{\,uv \;:\; u,v\in R(n)\ \text{disjoint as perfect sets}\,\}.
\]
\end{defn}

In this formulation a perfect partition is a clique of size $\frac{d_n}{n-1}$ in
$\Rcal_n$, and such a clique is necessarily of maximum size, since its members are
pairwise disjoint perfect sets of $n-1$ derangements each and so account for all derangements of order $n$. McKay compiled a list of all reduced Latin squares for
$n\le7$ \cite{M}, providing data for a computation search of all the  perfect sets, including the
elusive case $n=7$. From manual exhaustive search it is known that $L_{3,1}$ and
$L_{4,1}$ have exactly one perfect partition each \cite{BCL05}.

\subsection{The case \texorpdfstring{$n=5$}{n=5}}\label{subsec:n5}

We computed the maximum
cliques of $\Rcal_5$ by brute force. The search is exhaustive, and it returns
exactly $42$ perfect partitions of $L_{5,1}$. This is the first verification that
perfect partitions are non-unique: only one perfect partition of $L_{5,1}$ was
previously known \cite{BCL05}.

\begin{observation}\label{obs:pure}
Brualdi et al.\ found a perfect partition of $L_{5,1}$ consisting of $6$ pure
$(5)$-cycle perfect sets and $5$ pure $(3,2)$-cycle perfect sets \cite{BCL05}. Of the $42$ perfect partitions of $L_{5,1}$,
only $2$ keep the $(5)$ and $(3,2)$ cycle types in separate perfect sets. The other $40$ mix $(5)$ and $(3,2)$-cycles within a perfect set.
\end{observation}

\subsection{The case \texorpdfstring{$n=6$}{n=6}}\label{subsec:n6}

A greedy clique search algorithm on $\Rcal_6$ produced $82$ distinct perfect partitions of
$L_{6,1}$ before the search was manually terminated due to runtime, so it is likely that more exist. Previous work had exhibited a
single perfect partition of $L_{6,1}$ \cite{LSY12}. Together with the $n=5$ data, this supports non-uniqueness of perfect partitions, and the number of partitions appears to grow with $n$.

\subsection{The case \texorpdfstring{$n=7$}{n=7}}\label{subsec:n7}

The existence of perfect partitions for $n\ge7$ was left unresolved in
\cite{LSY12}. For $n=7$ the task is to partition the $1854$ derangements into
$309$ disjoint perfect sets of size $6$, drawn from a list of about
$16.9\times10^{6}$ reduced Latin squares.

A purely greedy search gets stuck at local optima. The one heuristic that made
substantial progress was a rarity-weighted greedy search, which incorporates the
rarest derangements first, covering  $97.4\%$ of $D_7$ but it did not complete a
partition. To obtain an exact cover we applied a hybrid of Knuth's Algorithm~X
(Dancing Links) \cite{K00} and the Bron--Kerbosch maximal clique algorithm
\cite{BK73} to the local neighbourhood of the highest-degree vertex of
$\Rcal_7$.

Two searches produced two perfect partitions.

\begin{enumerate}
  \item \textbf{$7$-cycle-free partition.} Since $7$ is prime,
        Lemma~\ref{lem:p-cycles} splits the $720$ $7$-cycles into $120$ disjoint
        perfect sets on their own. For this search we removed all reduced Latin squares containing a 7-cycle perfect matching. 
        The remaining $1134$ derangements were partitioned into $189$ perfect
        sets in our computational search, to form our complete perfect partition. We found some structure in the cycle-type combinations used to form perfect sets. In particular, $56$ of the sets
        are pure $(4,3)$, and $29$ sets consist of exactly two $(3,2,2)$
        derangements together with four $(5,2)$ derangements.
      \item \textbf{Mixed $7$-cycle partition.} This search looked through all reduced Latin squares for $n=7$, including those that use $7$-cycles. The $309$ disjoint perfect sets split into significantly more cycle types than the 7-cycle free case, forming $58$ distinct cycle-type
        profiles, of which only $6$ were pure. The mixing is thus much broader than the few types
        identified by Li et al. \cite{LSY12}, suggesting that heavily mixing cycle
        types within perfect sets is an effective mechanism for finding
        higher-order perfect partitions.
\end{enumerate}

We found no consistent algebraic action carrying one perfect set of these partitions to another. This lack of symmetry may explain the historical
difficulty of finding explicit algebraic constructions. From our computational search, despite the apparent abundance of perfect partitions, we did not identify a pattern for construction for a partition for an arbitrary $n$. 

\section{Future work}\label{sec:future}
From our computational results, we hypothesize partitions are generally not unique. 
\begin{conjecture}
    For all $n \geq 2$, there exists a perfect partition of $L_{n,1}$, and for all $n \geq 5$, there exist multiple perfect partitions of $L_{n,1}$.
\end{conjecture}

\subsection*{Improving the constructive bound}

The binding constraint in Section~\ref{sec:layers} is the block part, not the
layers. We illustrate this with $L_{10,1}$.

\begin{example}\label{ex:L10-gap}
Let $n=10$, so $p=5$, $k=2$, and there are $\tfrac{10!}{(5!)^2\,2!}=126$ vertex
partitions. For each partition, the layer $Q^1_{\V}$ supplies
$(5!)^2/5^{2}=576$ pairwise disjoint perfect sets by Theorem~\ref{thm:main}.
On the other hand, $L_{\V}$ is two copies of $L_{5,1}$, so it has $d_5^{\,2}=44^2$
perfect matchings and at most $44^2/4=484$ pairwise disjoint perfect sets.

Theorem~\ref{thm:general_lower_bound} uses $|\mathcal P(L_{5,1})|\cdot 3!\cdot 4=264$ perfect sets
of $L_{\V}$ for each partition, giving $126\cdot264=33\,264$ in total, or $22.4\%$
of $d_{10}/9=148\,329$. This can be improved by allowing our special set $\mathcal P(L_{p,1})$ to sit in
either block, rather than just $P_1$, while the other consists of only $p$-cycles. Let $\mathcal F$ be the perfect partition of $L_{5,1}$ with $6$
perfect sets made of $5$-cycles and $5$ made of $(3,2)$-cycles, and let
$P_1,P_2\in\mathcal F$ with at most one of them of type $(3,2)$. There are
$6^2+2\cdot5\cdot6=96$ such pairs, and with the $4$ rotations this gives $384$
perfect sets of $L_{\V}$, which is at most $576$. Two such sets are disjoint by
Lemmas~\ref{lem:diff-P} and~\ref{lem:diff-r}, since $\mathcal F$ is pairwise
disjoint. Sets from different partitions are disjoint because a matching with
cycle lengths $5,5$ or $5,3,2$ can be split into two groups of total length $5$
in only one way, and this split recovers $\V$. The total is
$126\cdot384=48\,384$ perfect sets, or $32.6\%$ of $d_{10}/9$.
\end{example}
The obstacle to going forward is that these sets can collide across partitions. A matching
whose cycles in the two blocks have types $(3,2)$ and $(3,2)$ has cycle lengths
$3,2,3,2$, which can be regrouped as $\{3,2\},\{3,2\}$ in a different way, so
the same matching also lies in $L_{\V'}$ for a second partition $\V'$. In the
construction above this never happens, because at most one block is not a
$5$-cycle. A natural next step is to decide, for each such matching, which of
the partitions containing it is allowed to keep it, while keeping every perfect
set whole.

A second route is to use more than one prime factor of $n$. For distinct prime
factors $p,p'$ of $n$, the block constructions of Section~\ref{sec:blocks} in
$\Lcal_{p,1;k}$ and $\Lcal_{p',1;k'}$ produce disjoint perfect sets. This is because every matching from the $p'$-construction
   contains a $p'$-cycle (in its second block), while every cycle of a matching
   from the $p$-construction has length at most $p<p'$.

The obstruction occurs at the $Q^\kappa_\V$ layer, where taut matchings of
the respective layer graphs are a shared resource when
$\gcd(k,\kappa)=\gcd(k',\kappa')$, so one would have to incorporate non-taut
matchings. Non-taut matchings, however, risk coinciding with each other or with matchings already used by the block construction, so such an extension requires care to avoid overlap.

\subsection*{Asymptotics via the hypergraph \texorpdfstring{$\Hcal_n$}{H(n)}}

A perfect partition of $L_{n,1}$ is a perfect matching of the $(n-1)$-uniform
hypergraph $\Hcal_n$ of Section~\ref{subsec:formulations}, which suggests
looking for one by nibble-type methods. The Pippenger--Spencer theorem
\cite{PS89} produces a matching covering all but $o(|V|)$ of the vertices of an
$r$-uniform hypergraph, provided $r$ is fixed, the degrees satisfy
$\delta/\Delta\to1$, and the maximum codegree is $o(\Delta)$. Three things stand in the way.

First, the uniformity of $\Hcal_n$ is $n-1$ and grows with $n$, so the theorem
does not apply to the family $(\Hcal_n)_{n\ge2}$ as it stands, and a version valid
for growing uniformity would be needed, and we are not aware of one.

Second, the degree condition is open. Part of it is provable: conjugation by
$g\in S_n$ permutes $D_n$ and maps perfect sets to perfect sets, so it is an
automorphism of $\Hcal_n$ and $\deg_{\Hcal_n}$ is constant on conjugacy
classes; only the $p(n)-p(n-1)$ cycle types need be compared. Moreover
$\deg_{\Hcal_n}(\sigma)=N(\sigma)/(n-2)!$, where $N(\sigma)$ is the number of
Latin squares with first two rows $\mathrm{id},\sigma$. Asymptotic equality of
these degrees is therefore a strengthening of the
conjecture of Cavenagh, Greenhill and Wanless~\cite{CGW08} on the distribution
of two rows of a random Latin square, and is open. Their result gives
agreement only up to a factor $n^{3/2}$ for almost all derangements.

Third, the codegree condition is open as well. If a near-perfect partition were obtained this way, the standard route to a full partition is with an absorbing structure \cite{RRS09}.

\section*{Partition Data}
The reduced Latin square partition data discussed in this paper is available on Zenodo: 

\href{https://doi.org/10.5281/zenodo.23503821}{https://doi.org/10.5281/zenodo.23503821}

\section*{Acknowledgments}
Computational resources were provided by William \& Mary Research Computing. The authors thank Chi-Kwong Li for helpful discussions. 

\section*{AI usage}

Opus 5.5 was used to proofread and revise the writing of the original manuscript, which was the second author's honor thesis (no AI was used). AI suggested to replace one $(p-2)!$ term with the term $|\mathcal P(L_{p,1})|$ (note that $|\mathcal P (L_{p,1})|\ge (p-2)!$) in the lower bound of Theorem~\ref{thm:general_lower_bound}, which we verified and adapted.

\section*{Appendix: \texorpdfstring{$7$}{7}-cycle-free perfect partition of \texorpdfstring{$L_{7,1}$}{L(7,1)}}\label{app:n7}

This appendix lists the $7$-cycle-free perfect partition of $L_{7,1}$ from
Section~\ref{subsec:n7}, consisting of $189$ pairwise disjoint perfect sets, partitioning all non-7-cycle derangements of order 7 into perfect sets. 

Note, we omitted $120$ perfect sets required to complete the full partition for $n=7$. The $720$ derangements of
cycle type $(7)$ are covered by the $120$ perfect sets of
Lemma~\ref{lem:p-cycles}, namely the sets $\langle\sigma\rangle\setminus\{e\}$
for each subgroup $\langle\sigma\rangle\le S_7$ of order $7$. These are fully
determined by the lemma, so we do not list them.

The remaining $1134$
derangements, none of which is a $7$-cycle, are partitioned into the $189$
perfect sets in Table~\ref{tab:rls7}. Each perfect set is given as a reduced
Latin square under the bijection of Lemma~\ref{lem:latin}. We also omit the identity element, $0123456$, in every reduced Latin square for compactness. The second perfect partition of $L_{7,1}$ described in Section~\ref{subsec:n7},
which uses $7$-cycles in mixed perfect sets, is too large to list here. It is
available, together with the partition below, in the data repository cited
under Partition Data.

{\scriptsize\setlength{\tabcolsep}{3pt}
\begin{longtable}{r l @{\hspace{14pt}} r l}
\caption{The $189$ perfect sets of the $7$-cycle-free perfect partition of
$L_{7,1}$ that contain no $7$-cycle, written as reduced Latin squares on the
symbols $0,\dots,6$. Each entry gives rows $1$--$6$; row $0$ is $0123456$ in
every square. The $120$ remaining perfect sets are those of
Lemma~\ref{lem:p-cycles}.}\label{tab:rls7}\\
\hline \# & Rows 1--6 & \# & Rows 1--6 \\ \hline \endfirsthead
\hline \# & Rows 1--6 & \# & Rows 1--6 \\ \hline \endhead
\hline \endfoot
1 & 1264530 2530164 3402615 4651203 5346021 6015342 & 2 & 1204635 2056314 3415260 4362501 5631042 6540123 \\ \hline
3 & 1452630 2015364 3264501 4306215 5641023 6530142 & 4 & 1254603 2540361 3462510 4601235 5316042 6035124 \\ \hline
5 & 1462503 2641035 3254610 4305261 5036142 6510324 & 6 & 1204563 2630145 3452601 4365210 5016324 6541032 \\ \hline
7 & 1564203 2016345 3402561 4251630 5630124 6345012 & 8 & 1542603 2461530 3254061 4036125 5610342 6305214 \\ \hline
9 & 1534260 2450631 3602145 4261503 5046312 6315024 & 10 & 1542360 2451603 3614025 4260531 5036214 6305142 \\ \hline
11 & 1540632 2451360 3264015 4612503 5306124 6035241 & 12 & 1504623 2456301 3065214 4231560 5640132 6312045 \\ \hline
13 & 1642530 2031645 3465201 4350162 5206314 6514023 & 14 & 1640523 2065341 3456210 4231605 5312064 6504132 \\ \hline
15 & 1452063 2541630 3615204 4036512 5260341 6304125 & 16 & 1465023 2536140 3604215 4250361 5041632 6312504 \\ \hline
17 & 1604532 2465310 3056241 4312605 5240163 6531024 & 18 & 1634520 2541063 3012645 4356201 5260134 6405312 \\ \hline
19 & 1340562 2435601 3061245 4256310 5604123 6512034 & 20 & 1340625 2436510 3051264 4265301 5614032 6502143 \\ \hline
21 & 1635024 2560143 3416205 4052631 5241360 6304512 & 22 & 1245063 2530641 3602514 4056132 5461320 6314205 \\ \hline
23 & 1364025 2516304 3651240 4035162 5240631 6402513 & 24 & 1302564 2461305 3654012 4035621 5246130 6510243 \\ \hline
25 & 1360524 2546031 3051642 4602315 5234160 6415203 & 26 & 1652034 2015643 3246510 4560321 5304162 6431205 \\ \hline
27 & 1206534 2651043 3045612 4530261 5364120 6412305 & 28 & 1435062 2561340 3652104 4016235 5204613 6340521 \\ \hline
29 & 1346520 2534061 3415602 4260315 5602134 6051243 & 30 & 1562304 2054613 3615042 4306521 5431260 6240135 \\ \hline
31 & 1230564 2645013 3516240 4062135 5304621 6451302 & 32 & 1536204 2615340 3451062 4360125 5042613 6204531 \\ \hline
33 & 1534602 2415063 3641520 4206135 5360214 6052341 & 34 & 1360245 2514603 3406512 4652130 5231064 6045321 \\ \hline
35 & 1246035 2460513 3605124 4051362 5312640 6534201 & 36 & 1560234 2345610 3054162 4236501 5612043 6401325 \\ \hline
37 & 1456032 2310645 3561204 4205163 5642310 6034521 & 38 & 1506324 2014635 3265041 4631502 5342160 6450213 \\ \hline
39 & 1562043 2430615 3245160 4316502 5601324 6054231 & 40 & 1460325 2641503 3514062 4235610 5306241 6052134 \\ \hline
41 & 1462035 2651304 3506142 4315620 5034261 6240513 & 42 & 1654230 2430561 3516024 4305612 5062143 6241305 \\ \hline
43 & 1605342 2046135 3564021 4352610 5410263 6231504 & 44 & 1634205 2561034 3416520 4052163 5340612 6205341 \\ \hline
45 & 1436025 2365140 3614502 4250613 5042361 6501234 & 46 & 1436502 2640315 3502641 4065123 5314260 6251034 \\ \hline
47 & 1435620 2016534 3504261 4650312 5261043 6342105 & 48 & 1645320 2416035 3201564 4562103 5034612 6350241 \\ \hline
49 & 1642305 2365014 3256140 4501263 5430621 6014532 & 50 & 1205643 2031564 3642015 4516320 5460231 6354102 \\ \hline
51 & 1405263 2614530 3250641 4536102 5361024 6042315 & 52 & 1302645 2650134 3041562 4516203 5264310 6435021 \\ \hline
53 & 1345602 2634015 3061524 4506231 5412360 6250143 & 54 & 1536042 2354160 3245601 4601523 5062314 6410235 \\ \hline
55 & 1350642 2514360 3461025 4605213 5236104 6042531 & 56 & 1652340 2345061 3214605 4506123 5061234 6430512 \\ \hline
57 & 1432605 2564130 3046521 4315062 5601243 6250314 & 58 & 1546302 2364510 3405621 4251063 5610234 6032145 \\ \hline
59 & 1065243 2354601 3502164 4610325 5436012 6241530 & 60 & 1230645 2014563 3645102 4561230 5406321 6352014 \\ \hline
61 & 1256304 2045163 3514620 4630215 5362041 6401532 & 62 & 1406235 2654103 3012564 4365021 5231640 6540312 \\ \hline
63 & 1564320 2630514 3201645 4015263 5346102 6452031 & 64 & 1265340 2356104 3542061 4610532 5401623 6034215 \\ \hline
65 & 1450362 2601534 3562140 4236015 5014623 6345201 & 66 & 1365204 2654310 3042165 4501632 5236041 6410523 \\ \hline
67 & 1650243 2346501 3064125 4512630 5401362 6235014 & 68 & 1405632 2310564 3546120 4652301 5064213 6231045 \\ \hline
69 & 1350264 2406315 3541602 4635120 5264031 6012543 & 70 & 1635240 2501364 3246105 4062513 5410632 6354021 \\ \hline
71 & 1036245 2341560 3560124 4615302 5402631 6254013 & 72 & 1356240 2634501 3041625 4510362 5462013 6205134 \\ \hline
73 & 1064523 2635104 3240615 4502361 5416230 6351042 & 74 & 1045623 2650341 3506214 4312560 5461032 6234105 \\ \hline
75 & 1502634 2356041 3410562 4061325 5634210 6245103 & 76 & 1354062 2061543 3406125 4615230 5230614 6542301 \\ \hline
77 & 1205364 2364105 3650241 4516032 5432610 6041523 & 78 & 1034562 2346015 3265104 4510623 5602341 6451230 \\ \hline
79 & 1042563 2614305 3256014 4531620 5360142 6405231 & 80 & 1342065 2064531 3650214 4215603 5436120 6501342 \\ \hline
81 & 1206345 2540613 3654201 4361520 5032164 6415032 & 82 & 1052643 2436105 3260514 4605321 5341062 6514230 \\ \hline
83 & 1235604 2041365 3560241 4652013 5406132 6314520 & 84 & 1234065 2560314 3605241 4351602 5046123 6412530 \\ \hline
85 & 1532064 2056143 3460215 4205631 5614320 6341502 & 86 & 1062534 2314065 3450612 4635201 5206143 6541320 \\ \hline
87 & 1406523 2351064 3610245 4562310 5234601 6045132 & 88 & 1245630 2064315 3650142 4512063 5436201 6301524 \\ \hline
89 & 1605234 2564013 3456102 4012365 5341620 6230541 & 90 & 1546023 2304561 3260145 4015632 5631204 6452310 \\ \hline
91 & 1054632 2401563 3562014 4216305 5630241 6345120 & 92 & 1465230 2506341 3052614 4210563 5634102 6341025 \\ \hline
93 & 1056342 2340165 3205614 4562031 5614203 6431520 & 94 & 1654302 2306145 3240561 4561023 5012634 6435210 \\ \hline
95 & 1062345 2415630 3540162 4206513 5634021 6351204 & 96 & 1056234 2361045 3214560 4502613 5640321 6435102 \\ \hline
97 & 1456203 2341065 3065142 4630521 5214630 6502314 & 98 & 1065324 2431065 3510642 4256103 5604231 6342510 \\ \hline
99 & 1506243 2341605 3215064 4632510 5064321 6450132 & 100 & 1046325 2305164 3652041 4560213 5431602 6214530 \\ \hline
101 & 1305624 2546103 3412065 4650231 5061342 6234510 & 102 & 1042635 2416503 3250164 4635012 5361240 6504321 \\ \hline
103 & 1035642 2506314 3450261 4612035 5264103 6341520 & 104 & 1034625 2450163 3546012 4602531 5261304 6315240 \\ \hline
105 & 1250634 2046513 3405162 4631025 5364201 6512340 & 106 & 1246503 2654031 3560142 4032615 5301264 6415320 \\ \hline
107 & 1064325 2351640 3510264 4605132 5246013 6432501 & 108 & 1362540 2536014 3450162 4605231 5241603 6014325 \\ \hline
109 & 1604325 2316540 3251604 4065231 5430162 6542013 & 110 & 1064235 2546310 3610524 4235061 5301642 6452103 \\ \hline
111 & 1402365 2635041 3560214 4016523 5341602 6254130 & 112 & 1046532 2410365 3501624 4265013 5634201 6352140 \\ \hline
113 & 1036524 2645301 3564210 4351062 5210643 6402135 & 114 & 1046523 2601345 3452160 4265031 5314602 6530214 \\ \hline
115 & 1256043 2564301 3640512 4302165 5031624 6415230 & 116 & 1035264 2456013 3610542 4301625 5264301 6542130 \\ \hline
117 & 1265034 2654301 3542610 4301562 5016243 6430125 & 118 & 1054263 2640531 3502614 4361025 5416302 6235140 \\ \hline
119 & 1052634 2435160 3260541 4316025 5641302 6504213 & 120 & 1065234 2456103 3241065 4532610 5604321 6310542 \\ \hline
121 & 1052364 2534610 3460521 4615032 5246103 6301245 & 122 & 1264305 2631540 3540621 4356012 5402163 6015234 \\ \hline
123 & 1645032 2036541 3512604 4351260 5460123 6204315 & 124 & 1035624 2360541 3204165 4651032 5416203 6542310 \\ \hline
125 & 1306542 2035614 3240165 4652031 5461203 6514320 & 126 & 1240365 2465031 3601542 4056213 5312604 6534120 \\ \hline
127 & 1465302 2016543 3640125 4352061 5201634 6534210 & 128 & 1042365 2610543 3504612 4256031 5361204 6435120 \\ \hline
129 & 1045362 2604531 3410625 4256013 5362104 6531240 & 130 & 1045632 2430165 3604521 4562013 5316204 6251340 \\ \hline
131 & 1204365 2065134 3540612 4632501 5416023 6351240 & 132 & 1630542 2546013 3015624 4361205 5204361 6452130 \\ \hline
133 & 1530624 2465103 3016542 4201365 5642031 6354210 & 134 & 1260543 2456130 3645021 4031265 5302614 6514302 \\ \hline
135 & 1430265 2516043 3642501 4265130 5301624 6054312 & 136 & 1056243 2564310 3640521 4312065 5231604 6405132 \\ \hline
137 & 1254360 2615034 3460125 4506213 5342601 6031542 & 138 & 1054362 2645031 3462105 4536210 5310624 6201543 \\ \hline
139 & 1546203 2014365 3250614 4635021 5462130 6301542 & 140 & 1205634 2361540 3410265 4056321 5642103 6534012 \\ \hline
141 & 1062543 2350614 3401265 4536021 5614302 6245130 & 142 & 1352604 2630541 3546210 4261035 5014362 6405123 \\ \hline
143 & 1602543 2054361 3465120 4536012 5210634 6341205 & 144 & 1304265 2051634 3460512 4615023 5246301 6532140 \\ \hline
145 & 1034265 2351604 3612540 4560312 5406123 6245031 & 146 & 1036542 2405163 3510624 4362015 5641230 6254301 \\ \hline
147 & 1654023 2510634 3062541 4301265 5436102 6245310 & 148 & 1206543 2645130 3014265 4530612 5462301 6351024 \\ \hline
149 & 1236540 2501634 3465012 4310265 5642301 6054123 & 150 & 1564302 2045631 3651024 4230165 5406213 6312540 \\ \hline
151 & 1340265 2015634 3654021 4506312 5462103 6231540 & 152 & 1650324 2034165 3206541 4365012 5412603 6541230 \\ \hline
153 & 1460532 2506143 3054621 4231065 5612304 6345210 & 154 & 1504362 2315604 3261540 4056123 5640231 6432015 \\ \hline
155 & 1640235 2301564 3465102 4516023 5032641 6254310 & 156 & 1056324 2460135 3645210 4531062 5214603 6302541 \\ \hline
157 & 1046235 2531604 3602541 4215063 5460312 6354120 & 158 & 1456023 2361504 3215640 4032165 5604312 6540231 \\ \hline
159 & 1645230 2530614 3402165 4216503 5361042 6054321 & 160 & 1032564 2465130 3210645 4651023 5346201 6504312 \\ \hline
161 & 1654320 2306514 3041265 4532601 5460132 6215043 & 162 & 1560342 2604513 3215604 4356021 5432160 6041235 \\ \hline
163 & 1432560 2605143 3546201 4210635 5364012 6051324 & 164 & 1532640 2604315 3210564 4365102 5046231 6451023 \\ \hline
165 & 1365042 2504163 3046215 4651320 5432601 6210534 & 166 & 1356024 2601543 3540261 4215630 5064312 6432105 \\ \hline
167 & 1346205 2605314 3451620 4560132 5214063 6032541 & 168 & 1632540 2405613 3056124 4561302 5340261 6214035 \\ \hline
169 & 1632504 2405631 3564012 4356120 5041263 6210345 & 170 & 1032645 2406513 3645201 4561032 5214360 6350124 \\ \hline
171 & 1456320 2501643 3640215 4032561 5364102 6215034 & 172 & 1354620 2506134 3641205 4065312 5432061 6210543 \\ \hline
173 & 1645203 2304165 3561042 4216530 5032614 6450321 & 174 & 1045263 2604135 3450621 4361502 5216340 6532014 \\ \hline
175 & 1564032 2405361 3216504 4051623 5632140 6340215 & 176 & 1532604 2401365 3654120 4360512 5216043 6045231 \\ \hline
177 & 1345260 2504631 3061542 4632015 5216304 6450123 & 178 & 1432065 2305614 3216540 4650321 5064132 6541203 \\ \hline
179 & 1054623 2605134 3461502 4210365 5346210 6532041 & 180 & 1546230 2305641 3214065 4650123 5461302 6032514 \\ \hline
181 & 1360542 2504613 3045261 4632105 5216034 6451320 & 182 & 1632045 2406531 3564102 4351620 5210364 6045213 \\ \hline
183 & 1065342 2406135 3541260 4350621 5632014 6214503 & 184 & 1540263 2304615 3456120 4061532 5632041 6215304 \\ \hline
185 & 1450623 2605341 3064512 4216035 5341260 6532104 & 186 & 1465032 2306541 3051624 4532160 5640213 6214305 \\ \hline
187 & 1350624 2504361 3216045 4065132 5641203 6432510 & 188 & 1364502 2401635 3650124 4532061 5046213 6215340 \\ \hline
189 & 1064532 2301645 3456021 4215360 5632104 6540213 &  &  \\ \hline
\end{longtable}}


\begin{thebibliography}{9}

\bibitem{BK73} C. Bron and J. Kerbosch, Algorithm 457: Finding all cliques of an
undirected graph, \emph{Commun. ACM} \textbf{16} (1973), 575--577.

\bibitem{BCL05} R. Brualdi, H. Chiang, C.-K. Li, A partition problem for sets of
permutation matrices, \emph{Bull. Inst. Combin. Appl.} \textbf{43} (2005),
67--79.

\bibitem{K00} D. E. Knuth, Dancing links, in \emph{Millennial Perspectives in
Computer Science}, Palgrave, 2000, pp.\ 187--214.

\bibitem{LSY12} C.-K. Li, J. Soosiah, G. Yu, Perfect partition of some regular
bipartite graphs, arXiv:1212.6091.

\bibitem{M} B. D. McKay, \emph{Combinatorial Data: Latin Squares},
\url{https://users.cecs.anu.edu.au/~bdm/data/latin.html}


\bibitem{PS89} N. Pippenger, J. Spencer, Asymptotic behavior of the chromatic
index for hypergraphs, \emph{J. Combin. Theory Ser. A} \textbf{51} (1989),
24--42.

\bibitem{CGW08} N. J. Cavenagh, C. Greenhill, I. M. Wanless, The cycle structure of
two rows in a random Latin square, \emph{Random Structures \& Algorithms}
\textbf{33} (2008), no.~3, 286--309.

\bibitem{RRS09} V. R\"odl, A. Ruci\'nski, and E. Szemer\'edi, Perfect matchings in
large uniform hypergraphs with large minimum collective degree, \emph{J. Combin.
Theory Ser. A} \textbf{116} (2009), no.~3, 613--636.

\end{thebibliography}
\end{document}